\documentclass{amsart}

\newtheorem{theorem}[equation]{Theorem}
\newtheorem{prop}[equation]{Proposition}
\newtheorem{lemma}[equation]{Lemma}

\theoremstyle{remark}
\newtheorem{remark}[equation]{Remark}
\theoremstyle{definition}

\numberwithin{equation}{subsection}

\renewcommand{\qed}{\hspace*{\fill} \setlength{\unitlength}{1mm}
\begin{picture}(2.5,2.5)
      \put(0,0){\framebox(2.5,2.5){}}
  \end{picture}
\setlength{\unitlength}{1pt}}

\newcommand{\Jac}{\operatorname{Jac}}
\newcommand{\supp}{\operatorname{Supp}}

\newcommand{\dist}{d_M}
\newcommand{\DIST}{d_{SM}}

\newcommand{\Id}{\operatorname{Id}}

\newcommand{\Horiz}{\operatorname{Horiz}}
\newcommand{\ver}{\operatorname{Vert}}

\newcommand{\integers}{{\bf Z}}

\newcommand{\reals}{{\bf R}}

\newcommand{\zed}{\integers}

\newcommand{\K}{{\mathcal{K}}}
\newcommand{\N}{{\mathcal{N}}}
\newcommand{\R}{{\mathcal{R}}}
\newcommand{\cH}{{\mathcal{H}}}
\begin{document}
\title[Estimates from below for the spectral function]
{Estimates from below for the spectral function and for the remainder
in local Weyl's law}
\author[D. Jakobson]{Dmitry Jakobson}
\address{Department of Mathematics and
Statistics, McGill University, 805 Sherbrooke Str. West,
Montr\'eal QC H3A 2K6, Ca\-na\-da.}
\email{jakobson@math.mcgill.ca}

\author[I. Polterovich]{Iosif Polterovich}
\address{D\'e\-par\-te\-ment de math\'ematiques et de
sta\-tistique, Univer\-sit\'e de Mont\-r\'eal CP 6128 succ
Centre-Ville, Mont\-r\'eal QC  H3C 3J7, Canada.}
\email{iossif@dms.umontreal.ca}


\thanks{The first author was supported by NSERC, FQRNT, Alfred P.
Sloan Foundation fellowship and Dawson fellowship.  The second
author was supported by NSERC and FQRNT}

\begin{abstract}
We obtain asymptotic lower bounds for the spectral function of the
Laplacian and for the remainder in local Weyl's law on
manifolds. In the negatively curved case,  thermodynamic formalism
is applied to improve the estimates. Key ingredients of the proof include the
wave equation parametrix, a pretrace formula and the Dirichlet box principle.
Our results develop and extend the unpublished thesis of A.~Karnaukh \cite{K}.
\end{abstract}
\maketitle


\section{Introduction and main results}
\subsection{Spectral function and Weyl's law}
Let $X$ be a compact Riemannian ma\-ni\-fold of dimension $n\ge 2$ with metric
$g_{ij}$ and of volume $V$. Let $\Delta$ be the Laplacian on
$X$ with the eigenvalues $0=\lambda_0<\lambda_1 \leq \lambda_2 \leq \dots$
and the corresponding orthonormal basis $\{\phi_i\}$ of eigenfunctions:
$\Delta \phi_i=\lambda_i\phi_i$.

Given $x,y\in X$, let
$$
N_{x,y}(\lambda)=\sum_{\sqrt{\lambda_i} \le \lambda}
\phi_i(x)\phi_i(y)
$$
be the spectral function of the Laplacian. On the  diagonal $x=y$
we denote it simply $N_x(\lambda)$. If
$N(\lambda)=\#\{\sqrt{\lambda_i} \le \lambda\}$ is the eigenvalue
counting function, then  $N(\lambda)=\int_X N_x(\lambda) dV$.
Let
\begin{equation}
\label{disk}
\sigma_n=\frac{2\pi^{n/2}}{n\Gamma(n/2)}
\end{equation}
be the volume of the unit ball in $\reals^n$. The asymptotic
behavior of the spectral and the counting functions is given by
(\cite{H}, see also \cite{Shubin}):
\begin{equation}\label{Weyl}\begin{aligned}
N_{x,y}(\lambda)&=O(\lambda^{n-1}),\qquad x\neq y;\\
N_x(\lambda)&=\frac{\sigma_n}{(2\pi)^n}\lambda^n+
R_x(\lambda),\qquad R_x(\lambda)=O(\lambda^{n-1});\\
N(\lambda)&=\frac{V\sigma_n}{(2\pi)^n}\lambda^n + R(\lambda),
\qquad R(\lambda)=O(\lambda^{n-1}).
\end{aligned}
\end{equation}
We refer to the asymptotics of $N_x(\lambda)$ as {\it local Weyl's
law}, the asymptotics for $N(\lambda)$ being the usual Weyl's law
for the distribution of eigenvalues.

In the present paper we focus on asymptotic  {\it lower} bounds
for the spectral function and for the remainder in local Weyl's law.
We recall that $f_1(\lambda)=\Omega(f_2(\lambda))$  for a
function $f_1$ and a positive function $f_2$ means
$\limsup_{\lambda \to \infty} |f_1(\lambda)|/f_2(\lambda)>0$.
\begin{theorem}
\label{weak} Let $X$ be a compact $n$-dimensional Riemannian
manifold, and let $x,y \in X$  be two points that are not
conjugate along any shortest geodesic joining them. Then
\begin{equation}\label{bound:weak}
N_{x,y}(\lambda)=\Omega\left(\lambda^{\frac{n-1}{2}}\right).
\end{equation}
\end{theorem}

Let us now formulate the on-diagonal counterpart of Theorem
\ref{weak}. Consider the heat trace asymptotics as $t \to 0^+$:
\begin{equation}
\label{heat} \sum_i e^{-\lambda_i t} \sim \frac{1}{(4\pi)^{n/2}}
\sum_{j=0}^\infty \int_X a_j(x) dV \,\,\, t^{j-\frac{n}{2}},
\end{equation}
where the local heat invariants $a_j(x)$ are the coefficients in
the asymptotic expansion of the heat kernel (see section
\ref{sec:new}). We shall use the following notation. Let
$\kappa_x=\min\{j\ge 1|\,\, a_j(x) \neq 0\}$. If $a_j(x)=0$ for
all $j\ge 1$ we set $\kappa_x=\infty$. Note that
$a_1(x)=\frac{\tau(x)}{6}$, where $\tau(x)$ is the scalar
curvature of $X$ at the point $x$.
\begin{theorem}
\label{weakon} Let $X$ be an $n$-dimensional Riemannian manifold and
$x \in X$ be an arbitrary point.

(i) If $n-2\kappa_x > 0$ then
\begin{equation}
\label{bound:weakon} R_x(\lambda)=\Omega(\lambda^{n-2\kappa_x}).
\end{equation}

(ii) If $n-4\kappa_x+1 < 0$ and if $x$ is not conjugate to itself
along any shortest geodesic loop, then
\begin{equation}
\label{weakncp} R_x(\lambda)=\Omega(\lambda^{\frac{n-1}{2}}).
\end{equation}
\end{theorem}
\begin{remark}
If the scalar curvature $\tau(x) \neq 0$ then $\kappa_x=1$ and
\eqref{bound:weakon} becomes $R_x(\lambda)=\Omega(\lambda^{n-2})$.
If $n-4\kappa_x+1 < 0$  then \eqref{weakncp} gives a stronger
estimate than \eqref{bound:weakon}. A bound similar to
\eqref{bound:weakon} holds also for the integrated remainder
$R(\lambda)$ in the usual Weyl's law, see section \ref{sec:new}.
\end{remark}
Theorems \ref{weak} and \ref{weakon} (ii) are proved in sections
\ref{35} and \ref{36} by analyzing the asymptotics of the parametrix
for the wave equation. The proof of Theorem \ref{weakon} (i) in
section \ref{21}  is based on the heat kernel asymptotics. Estimate
\eqref{weakncp} in dimension two should be compared with the
classical Hardy-Landau bound (and its generalization in
\cite[Proposition 3.1]{Sar}) $R(\lambda)=\Omega(\sqrt{\lambda})$ for
the remainder in the Gauss circle problem or, equivalently, for the
remainder in Weyl's law on a $2$-dimensional flat square torus (see
\cite{Sound} for the most recent improvement of such a bound). Note
that on a torus $R_x(\lambda) \equiv R(\lambda)$.

\subsection{Oscillatory error term in Weyl's law}
Let us represent the local Weyl's law in a somewhat different
form, introducing the {\it oscillatory error term}
$R_x^{osc}(\lambda)$:
\begin{equation}
\label{osc} N_x(\lambda)=\frac{1}{(4\pi)^{n/2}}
\sum_{j=0}^{\left[\frac{n-1}{2}\right]}
\frac{a_j(x)}{\Gamma\left(\frac{n}{2}-j+1\right)}\lambda^{n-2j}+R^{osc}_x(\lambda).
\end{equation}
This is {\it not} an asymptotic expansion since apriori
$R_x^{osc}(\lambda)=O(\lambda^{n-1})$ as the usual remainder
$R_x(\lambda)$.  The oscillatory remainder satisfies certain
``logarithmic Gaussian error estimates'', introduced many years ago
in \cite{Brown} . This approach was adopted by physicists and the
representation \eqref{osc} of Weyl's law often appears in physics
literature ( cf. \cite[(IV.13), p.37]{BH}, see also \cite{BB}). The
idea is to subtract the contribution of the ``singularity at zero''
(that we took into account in \eqref{bound:weakon}) to the counting
function.

\begin{theorem}\label{thm:osc}
Let $X$ be an $n$-dimensional Riemannian manifold and let $x \in
X$ be not conjugate to itself along any shortest geodesic loop.
Then
$$R_x^{osc}(\lambda) = \Omega(\lambda^{\frac{n-1}{2}}).$$
\end{theorem}

Theorem \ref{thm:osc} in proved in section \ref{pf:osc}.
\subsection{Thermodynamic formalism}
Asymptotic lower bounds \eqref{bound:weak} \eqref{bound:weakon}
and \eqref{weakncp} can be improved for manifolds of negative
curvature, see section \ref{neg}. We assume that for any pair of
directions $\xi, \eta$ the sectional curvature $K(\xi, \eta)$
satisfies
\begin{equation}
\label{curvature}
-K_1^2 \le K(\xi, \eta) \le -K_2^2.
\end{equation}

Apart from the usual wave equation techniques (cf. \cite{DG},
\cite{Berard}, \cite{K}) our
method uses the thermodynamic formalism (see, for example, \cite{Bowen},
\cite{PP}).
Let $G^t$ be the geodesic flow on the unit tangent bundle $SX$ and let
$E_\xi^u$ be the unstable subspace for $G^t$, $\xi \in SX$. The
{\it Sinai-Ruelle-Bowen potential} is a H\"older continuous
function $\cH:SX \to\reals$ which for any $\xi\in SX$ is defined
by the formula (see \cite{BR}, \cite{Sinai2})
\begin{equation}\label{SRB}
\cH(\xi)=\left.\frac{d}{dt}\right|_{t=0}\ln\det dG^t|_{E_\xi^u},
\end{equation}
For any continuous function $f:SX \to \reals$ one can define the
{\it topological pressure}
\begin{equation}\label{pressure}
P(f)=\sup_\mu \left(h_\mu + \int f d\mu \right),
\end{equation}
where the supremum is taken over all $G^t$-invariant measures $\mu$
and $h_\mu$ denotes the {\it measure-theoretic entropy} of the
geodesic flow (see \cite{Bowen}).
 In particular $P(0)=h$,
where $h$ is the {\it topological entropy} of the flow. It is
well-known that for the Sinai-Ruelle-Bowen potential $P(-\cH)=0$
and the corresponding equilibrium measure (i.e. the measure on
which the supremum is attained) is the Liouville measure $\mu_L$
on the unit tangent bundle:
\begin{equation}\label{ent}
h_{\mu_L} = \int_{SX} \cH d\mu_L
\end{equation}


\subsection{Estimates for negatively curved manifolds} \label{neg}
We now present our main results for manifolds of negative curvature.
\begin{theorem}\label{main:offdiag}
On a compact $n$-dimensional negatively curved manifold the spectral
function $N_{x,y}(\lambda)$ satisfies  for any $\delta > 0$ and
$x\neq y$:
\begin{equation}\label{bound:offdiag}
N_{x,y}(\lambda)=\Omega\left(\lambda^{\frac{n-1}{2}}\,
(\log\lambda)^{\frac{P(-\cH/2)}{h}-\delta}\right).
\end{equation}
\end{theorem}

\begin{theorem}\label{main}
The remainder in the local Weyl's law on an $n$-dimensional
compact negatively curved manifold satisfies:
\begin{equation}
\label{bound} R_x(\lambda)=\begin{cases} \Omega
\left(\lambda^{\frac{n-1}{2}}\,
(\log\lambda)^{\frac{P(-\cH/2)}{h}-\delta}\right) \,\,\, \forall \,
\delta >0, \quad n=2, 3\,;\cr \Omega (\lambda^{n-2}), \quad n \ge
4.\cr
\end{cases}
\end{equation} Also, in any dimension
\begin{equation} \label{oscneg} R_x^{osc}(\lambda)=\Omega
\left(\lambda^{\frac{n-1}{2}}\,
(\log\lambda)^{\frac{P(-\cH/2)}{h}-\delta}\right) \,\,\, \forall
\, \delta >0.
\end{equation}
\end{theorem}
\noindent For $n=2$, a bound similar to \eqref{bound} is proved in
\cite{K}, see a discussion below.

Theorems \ref{main:offdiag} and  \ref{main} are proved in section
\ref{proofs}.

Let us estimate the power $\frac{P(-\cH/2)}{h}$ in terms of
curvatures.  It was shown in \cite{Sinai1} that on an
$n$-dimensional negatively curved manifold satisfying
\eqref{curvature},
$$
(n-1)K_2\le h_{\mu_L} \le h \le (n-1)K_1
$$
we can obtain (using the definition of pressure and \eqref{ent}):
\begin{equation}\label{press:est}
P\left(-\frac{\cH}{2}\right)\ge h_{\mu_L} - \int \cH \mu_L +
\frac{1}{2} \int \cH d\mu_L=\frac{h_{\mu_L}}{2} \ge
\frac{(n-1)K_2}{2}
\end{equation}
and therefore
$$
\frac{P\left(-\cH/2\right)}{h} \ge \frac{K_2}{2K_1}>0.
$$

\subsection{Discussion}
\label{discussion}
Remainder estimates  in \eqref{Weyl} for $N(\lambda)$ and  $N_x(\lambda)$
are attained for round
spheres and hence are sharp. One can replace $O(\lambda^{n-1})$ by
$o(\lambda^{n-1})$
in the bounds for $R(\lambda)$ (respectively, for $R_x(\lambda)$) if the
initial directions
of closed geodesics (respectively, of geodesic loops at $x$) form a set of
measure zero
(for $R(\lambda)$ see \cite{DG}, for $R_x(\lambda)$ see \cite{Saf},
\cite[section 1.8]{SV}, and also \cite{SZ}).

Both local and integrated remainder estimates for Weyl's law
were actively studied under various geometric conditions (see \cite{CdV},
\cite{Ivrii}, \cite{Vol},
\cite{TP}, etc.).
Of particular interest are manifolds of negative (or, more generally, nonpositive)
curvature.
In this case,  $R_x(\lambda)=O(\lambda^{n-1}/\log \lambda)$ (\cite{Berard}).
It is conjectured in \cite{Randol} that on a surface of constant negative
curvature
$R(\lambda)=O(\lambda^{\frac{1}{2}+\epsilon})$ for any $\epsilon>0$.
In the same article it is proved that an average over $\lambda$ of
$|N_{x,y}(\lambda)|$ for
$x\neq y$ satisfies such an upper bound. On surfaces corresponding
to quaternionic lattices $R(\lambda)=\Omega
\left(\frac{\sqrt{\lambda}}{\log \lambda}\right)$
(see \cite{Hejhal}). Hence, if there are no additional assumptions, the bound
conjectured in
\cite{Randol} is almost the best possible.

The present paper originated as an attempt to improve the estimate
$R_x(\lambda)=\Omega(\sqrt{\lambda})$
obtained in  an unpublished Princeton Ph.D. thesis \cite{K} for
surfaces of arbitrary negative curvature.
The main novelty of our approach compared to \cite{K} is the use of
thermodynamic formalism.
It is also noticed that the restriction to
the negatively curved case considered in \cite{K} can be replaced
by much more general assumptions, see Theorems \ref{weak} and
\ref{weakon}.  We estimate
$N_{x,y}(\lambda)$ not just on the diagonal, but also for $x\neq y$.
Our results apply to arbitrary dimensions, while \cite{K} deals only
with surfaces.
Working in higher dimensions makes the analysis of the parametrix more
difficult.
In particular, asymptotics have different nature in high and low
dimensions, see Theorem \ref{main}.

On surfaces of constant negative curvature $\eqref{bound}$
coincides with the bound in \cite{K}, but the techniques
of thermodynamic formalism for hyperbolic flows allow us to
improve the results obtained in \cite{K} for manifolds of variable
negative curvature.
We do not need the hypothesis $K_1/K_2<2$ of \cite{K} to get a
logarithmic improvement in the estimates and, moreover, we get higher powers of
the logarithm. Note that the results of \cite{Randol} cited above make it
unlikely that in
two dimensions the power of $\lambda$ in a lower bound for $N_x(\lambda)$
can be made
greater than $1/2$. We would like to emphasize that the techniques of hyperbolic dynamics
can be quite useful
for proving spectral estimates on negatively curved manifolds  (see also
\cite{Vol} where results of
Bowen were applied to prove upper bounds for the remainder in Weyl's law).

Let us conclude by mentioning that on a {\it generic} negatively curved
surface it is believed that $R(\lambda)=O(\lambda^{\epsilon})$
for any $\epsilon>0$.  Such an
estimate looks plausible in view of the results on spectral
fluctuations, e.g. in \cite{Berry}, \cite{BS} and \cite{ABS}. We
note the difference between the predicted upper bound for the {\em
global} error term, and the lower bounds in this paper for
the {\em local} remainder.

Main results of the present paper were announced in \cite{JP}.

\section{Heat equation and wave equation techniques}
\label{sec:pretrace}
\subsection{A heat kernel approach to lower bounds}
\label{sec:new} In this section we prove part (i) of the Theorem
\ref{weakon}. While the rest of the paper is based on the wave
equation techniques, the contribution of the ``singularity at zero''
to the remainder in Weyl's law can be easier seen through the heat
asymptotics.

\smallskip

\noindent {\bf Proof of Theorem \ref{weakon} (i).}
\label{21}
For simplicity we
shall assume that $\kappa_x=1$, if not the proof is analogous.
Consider the short time asymptotic expansion of the heat kernel on
the diagonal
\begin{equation} \label{heatexp} \K(t,x,x)=\sum_i
e^{-\lambda_i t} \phi_i^2(x) \sim_{t\to 0^+}
\frac{1}{(4\pi)^{n/2}} \sum_{j=0}^{\infty} a_j(x)
t^{j-\frac{n}{2}}.
\end{equation}
The first terms are $a_0(x)=1$, $a_1(x)=\tau(x)/6$, $\tau(x)$ is
the scalar curvature at $x$. As we assumed above, $a_1(x) \neq 0$.

To simplify calculations we use the renormalized counting function
\begin{equation}
\label{renorm} \N_x(\mu)=\sum_{\lambda_i \le \mu} \phi_i^2(x)
=\frac{\sigma_n}{(2\pi)^n} \mu^{n/2} + \R_x(\mu),
\end{equation}
where $\R_x(\mu)$ is the renormalized remainder.
 Rewriting $\K(t,x,x)=\int_0^{\infty} e^{-t\mu}
d\N_x(\mu)$ and integrating by parts we get:
$$\K(t,x,x)=t\int_0^{\infty} e^{-t\mu} \N_x(\mu)
d\mu=t\int_0^{\infty}\frac{\sigma_n}{(2\pi)^n}
e^{-t\mu}\mu^{n/2}d\mu + t\int_0^{\infty}e^{-t\mu} \R_x(\mu)d\mu,$$
Comparing this with \eqref{heatexp} and taking into account
\eqref{disk} we note that the contribution of the principal term
cancels out the contribution of the $0$-th heat invariant to
$\K(t,x,x)$. Therefore,
\begin{equation}
\label{ostatok} \int_0^{\infty} e^{-t\mu} \R_x(\mu) d\mu =
a_1(x)t^{-\frac{n}{2}} + O\left(t^{1-\frac{n}{2}}\right), \quad
t\to 0^+.
\end{equation}
We shall prove by contradiction that \eqref{ostatok} implies
$\R_x(\mu)=\Omega(\mu^{\frac{n}{2}-1})$. Indeed, suppose that for
any $\varepsilon
> 0$ there exists $z=z(\varepsilon)$ such that
$|\R_x(\mu)|<\varepsilon \mu^{\frac{n}{2}-1}$ if $\mu \ge z$. We
have:
$$
\left|\int_0^{\infty} e^{-t\mu} \R_x(\mu)d\mu\right| \le
\left|\int_0^z e^{-t\mu} \R_x(\mu)d\mu \right| +
\left|\int_z^{\infty} e^{-t\mu} \R_x(\mu)d\mu \right|= |I_1|
+|I_2|.$$ The first term $|I_1| < const$, where the constant
depends on $z$ but not on $t$.  The second term can be estimated
by
$$
|I_2| \le \int_z^{\infty} \varepsilon e^{-t\mu}
\mu^{\frac{n}{2}-1}d\mu = \varepsilon
t^{-\frac{n}{2}}\int_{zt}^{\infty} e^{-\nu}
\nu^{\frac{n}{2}-1}d\nu\le \varepsilon
t^{-\frac{n}{2}}\int_{0}^{\infty} e^{-\nu} \nu^{\frac{n}{2}-1}d\nu =
\varepsilon t^{-\frac{n}{2}}\Gamma\left(\frac{n}{2}\right).
$$
Choosing $\epsilon < \frac{|a_1(x)|}{2 \Gamma(n/2)}$ we get
$$\left|\int_0^{\infty} e^{-t\mu} \R_x(\mu) d\mu \right| \le
\frac{|a_1(x)|}{2} t^{-\frac{n}{2}}$$
as $t \to 0^+$, which contradicts \eqref{ostatok}. Therefore,
$\R_x(\mu)=\Omega(\mu^{\frac{n}{2}-1})$. Comparing \eqref{renorm}
and \eqref{Weyl} we see that it is  equivalent to
$R_x(\lambda)=\Omega(\lambda^{n-2})$. This completes the proof of
the theorem. \qed

In order to prove a similar bound for the integrated remainder
$R(\lambda)$ one should just repeat the same argument for
$N(\lambda)$ instead of $N_x(\lambda)$ and for the heat trace
instead of the heat kernel.
\begin{remark} In the earlier version of the paper we proved a
weaker bound $R_x(\lambda)=\Omega(\lambda^{n-2\kappa_x-1})$ using
Proposition \ref{leadingon}. It was pointed out to us by Yu.
Safarov that a better estimate \eqref{bound:weakon} should be
valid. He also suggested a somewhat different proof of this result
based on the asymptotics of the Riesz means, cf. \cite{Saf2}.
\end{remark}

\subsection{Smoothed Fourier transform of the wave kernel}
The even part of the wave
kernel on $X$ satisfies
\begin{equation}\label{waveker1}
e(t,x,y)\ =\ \sum_{i=0}^\infty \cos(\sqrt{\lambda_i}t)
\phi_i(x)\phi_i(y)
\end{equation}
It is the fundamental solution of the wave equation
\begin{equation}
\label{waveqn}
\begin{aligned}
(\partial^2/\partial t^2-\Delta)e(t,x,y)&=0,\\
e(0,x,y)&=\delta(x-y),\\
(\partial/\partial t)e(0,x,y)&=0.
\end{aligned}
\end{equation}

Take a smooth function $\psi\in C_0^\infty(\reals)$ such that ${\rm
supp}\;\psi\subseteq [-1,1]$, it is  even and monotone decreasing on
[0,1], $\psi\geq 0$, $\psi(0)=1.$ Fix two large positive parameters
$\lambda,T$ and choose a compactly supported smooth test function
\begin{equation}\label{testfxn}
(1/T)\psi(t/T)\cos(\lambda t).
\end{equation}
Given $x,y\in M$, denote the integral in the right-hand side by
\begin{equation}\label{klamt}
k_{\lambda,T}(x,y)=\int_{-\infty}^\infty\frac{\psi(t/T)}{T}
\cos(\lambda t) e(t,x,y) dt
\end{equation}
Substituting \eqref{waveker1} into \eqref{klamt} we obtain
\begin{equation}\label{KH}
k_{\lambda, T}(x,y)=\sum_{i=0}^\infty \phi_i(x)\phi_i(y)
H_{\lambda,T}(\sqrt{\lambda_i}),
\end{equation}
where
\begin{equation}\label{defH}
H_{\lambda,T}(r)=\int_{-\infty}^\infty\frac{\psi(t/T)}{T}
\cos(\lambda t) \cos(r t)
dt=\frac{\widehat{\psi}(T(r+\lambda))+\widehat{\psi}(T(r-\lambda))}{2}.
\end{equation}
Here
$$
\widehat{\psi}(s)=\int_{-\infty}^{\infty}
e^{is\zeta}\psi(\zeta)d\zeta
$$
is the Fourier transform of $\psi$. Replacing the sum in
\eqref{KH} by an integral, we get the following representation of
$k_{\lambda,T}(x,y)$:
\begin{multline}\label{klamt2}
 \int_0^\infty H_{\lambda,T}(r)
dN_{x,y}(r)=\\
\int_0^\infty
\frac{\widehat{\psi}(T(r+\lambda))+\widehat{\psi}(T(r-\lambda))}{2}dN_{x,y}(r)=
k_{\lambda,T}(x, y)
\end{multline}
Formula \eqref{klamt2} plays a key role in our analysis.

On the diagonal, we shall use the following notation:
\begin{equation}\label{klamt3}
\int_0^\infty H_{\lambda,T}(r)
dR_x(r)= \widetilde{k}_{\lambda,T}(x)
\end{equation}
Note that the contribution of the main term in Weyl's law has been subtracted.

\subsection{Relation between $N_{x,y}(\lambda)$ and $k_{\lambda, T}(x,y)$,
$R_x(\lambda)$ and $\widetilde{k}_{\lambda,T}(x)$}
Assume that $T$ is either fixed or $T(\lambda) \to \infty$ as
$\lambda \to \infty$. The following lemma will be used to prove
Theorems \ref{weak} and \ref{weakon}.
\label{int:above}
\begin{lemma}
\label{above:weak1}
(i) Let $N_{x,y}(\lambda)=o(\lambda^a)$, $a>0$. Then
$k_{\lambda, T}(x,y)=o(\lambda^a)$. (ii) Let $R_x(\lambda)=o(\lambda^a)$.
Then $\widetilde{k}_{\lambda,T}(x)=o(\lambda^a)$.
\end{lemma}
\noindent{\bf Proof.} We shall prove (i), the proof of part (ii) is analogous.
Denote $\widehat{\psi}\,'$ by $\Psi$.
Since $\psi$ has compact support, $\Psi$ is Schwartz class.
By the assumption of the lemma, for any $\epsilon>0$,
$N_{x,y}(\lambda)< \epsilon\lambda^a$ for large enough $\lambda$.
Consider the left hand side of \eqref{klamt2}:
\begin{equation}
\label{LL}
\int_0^\infty H_{\lambda,T}(r) dN_{x,y}(r).
\end{equation}
Taking into account \eqref{klamt2} and  integrating \eqref{LL} by parts
we obtain
\begin{equation}
\label{byparts}
k_{\lambda, T}(x,y)\le \frac{\epsilon T}{2}\int_0^\infty |\Psi(T(r-\lambda))|
r^a dr +\frac{\epsilon T}{2}\int_0^\infty |\Psi(T(r+\lambda))|
r^a dr.
\end{equation}
Since $\Psi$ is Schwartz class, the second term of \eqref{byparts}
is $O(1)$.  Changing variables in the first term of
\eqref{byparts},
we obtain
$$
\frac{\epsilon T}{2}\int_{-\lambda}^{\infty} |\Psi(Ts)| (s+\lambda)^a ds
= \frac{\epsilon \lambda^a}{2}\int_{-\lambda T}^{\infty} |\Psi(u)|
\left(1+\frac{u}{\lambda T}\right)^a du \le C \epsilon \lambda^a
$$
for some constant $C>0$, where the last inequality again follows from
the fact that $\Psi$ is Schwartz class. Since $\epsilon$ can be chosen
arbitrarily small, we get $k_{\lambda, T}(x,y)=o(\lambda^a)$,
and this completes the proof of the lemma. \qed

\smallskip

A modification of this lemma is used to prove Theorems
\ref{main} and \ref{main:offdiag}.
\begin{lemma}
\label{above2}
(i) Let $N_{x,y}(\lambda)=O(\lambda^a (\log \lambda)^b)$, $a,b>0$. Then
$k_{\lambda, T}(x,y)=O(\lambda^a (\log \lambda)^b)$.
(ii) Let $R_x(\lambda)=O(\lambda^a (\log \lambda)^b)$, $a,b>0$. Then
$\widetilde{k}_{\lambda, T}(x)=O(\lambda^a (\log \lambda)^b)$.
\end{lemma}
\noindent{\bf Proof.}
This lemma is proved similarly to the previous
one. Again, we shall prove (i), the proof of (ii) is analogous.

By the assumption of the lemma there exists a constant $C>0$ such
that
$N_{x,y}(\lambda)<C \lambda^a (\log(1+\lambda))^b$. Writing a
representation of $k_{\lambda,T}(x,y)$ analogous to \eqref{byparts} we
get that the second term is $O(1)$ as before, and the first term can
be rewritten as
\begin{multline*}
\frac{C T}{2}\int_{-\lambda}^{\infty} |\Psi(Ts)| (s+\lambda)^a
(\log(1+\lambda+s))^b ds=\\\frac{C \lambda^a (\log\lambda)^b}{2}
\int_{-\lambda T}^{\infty} |\Psi(u)|
\left(1+\frac{u}{\lambda T}\right)^a
\left(1+\frac{\log(1+\frac{1+u/T}{\lambda})}{\log\lambda}\right)^bdu=
O(\lambda^a (\log\lambda)^b).
\end{multline*}
\begin{remark} Statements similar to
Lemmas \ref{above:weak1} and \ref{above2} are proved in \cite{K},
see also \cite[p.~226]{Sar}.
\end{remark}


\subsection{The pretrace formula}
If $X$ is negatively curved it is useful to consider the
fundamental solution  $E(t,x,y)$ of the wave equation
\eqref{waveqn} on the universal cover $M$ of $X$. Then given
$x,y\in X$, we have
\begin{equation}
\label{wavesum1} e(t,x,y)=\sum_{\omega\in\Gamma}E(t,x, \omega y),
\end{equation}
where the sum is taken over $\Gamma=\pi_1(X)$.
Let $K_{\lambda,
T}(x,y)$ be the analogue of $k_{\lambda,T}(x,y)$ corresponding to the wave kernel
$E(t,x,y)$ on $M$:
\begin{equation}
\label{klamt1}
K_{\lambda,T}(x,y)=\int_{-\infty}^\infty\frac{\psi(t/T)}{T}
\cos(\lambda t) E(t,x,y) dt
\end{equation}

Then \eqref{klamt2} becomes
\begin{equation}\label{pretrace2}
k_{\lambda,T}(x,y)=\int_0^\infty H_{\lambda,T}(r)
dN_{x,y}(r)=
\sum_{\omega\in\Gamma}K_{\lambda,T}(x,\omega y).
\end{equation}

\section{Asymptotics of the smoothed Fourier transform}
\subsection{Parametrix}\label{parametrix}
In this section we review the Hadamard parametrix construction for
the wave equation  (cf. \cite[Prop. 27]{Berard}, see also
\cite{Zelditch}). We present it in the setting when $X$ is a
manifold without conjugate points (as in the negatively curved
case) and hence the parametrix is globally defined on the
universal cover $M$.

 We shall work in
Riemannian normal coordinates centered at $x\in M$. Given $x,y\in
M$, let $r=d(x,y)$ and let $E(t,x,y)$ be the fundamental solution
of \eqref{waveqn}  on $M$. The parametrix for $E(t,x,y)$ is given
by:
\begin{equation}\label{param}
E(t,x,y)=\frac{1}{\pi^{\frac{n-1}{2}}} |t|
\sum_{j=0}^{\infty} u_j(x,y)
\frac{(r^2-t^2)_-^{j-\frac{n-3}{2}-2}} {4^j
\Gamma(j-\frac{n-3}{2}-1)} \mod C^{\infty},
\end{equation}
The expression (\ref{param}) is understood in the sense of
generalized  functions \cite{Gelfand}. The coefficients $u_j(x,y)$
are the solutions of the transport equations (\cite{Berard}) along
the geodesic joining $x$ and $y$ (since $X$ is negatively curved,
such a geodesic on $M$ is unique). In particular
$u_0(x,y)=g^{-1/2}(y)$, where $g=\sqrt{\det g_{ij} }$. The
on-diagonal values $u_j(x,x)=a_j(x)$ are the local heat invariants,
see section \ref{sec:new}.

We substitute \eqref{param} into \eqref{klamt} to get
\begin{equation}\label{klamt:param}\begin{aligned}
K_{\lambda,T}(x,y)&=
\int_{-\infty}^\infty\frac{\psi(t/T)
\cos(\lambda t)|t|}{T\pi^{\frac{n-1}{2}}}
\sum_{j=0}^{J}u_j(x,y) \frac{(r^2-t^2)_-^{j-\frac{n-3}{2}-2}}
{4^j\Gamma(j-\frac{n-3}{2}-1)}dt\\\;&+E_{\lambda,T}(x,y)
\end{aligned}\end{equation}
where  $J\ge 2([n/2]+1)+1$ and
$E_{\lambda,T}(x,y)=O(1)\exp(O(T)).$  (\cite[p.266]{Berard}).


\subsection{Leading term cancellation in the pretrace formula}
The results of sections 3.2--3.4 are essentially not new and could
be deduced from \cite{Berard}. We refer to \cite[p. 46]{DG} and
\cite[p. 259]{Berard} for the proper regularization of $E(t,x,x)$ at
$t=0$.
\begin{lemma}
\label{diag:claim}
The following relation holds:
\begin{multline}
\label{pretrem} \int_0^\infty H_{\lambda,T}(r) dR_x(r)=
\int_{-\infty}^\infty\frac{\psi(t/T) \cos(\lambda
t)}{T\pi^{\frac{n-1}{2}}}
\sum_{j=1}^{J}\frac{u_j(x,x)|t|^{2j-n}}{4^j
\Gamma(j-\frac{n-3}{2}-1)}dt +\\
\sum_{g\in \Gamma\setminus\{\operatorname{Id}\}}
\int_{-\infty}^\infty\frac{\psi(t/T) \cos(\lambda
t)|t|}{T\pi^{\frac{n-1}{2}}} \sum_{j=0}^{J}\frac{u_j(x,\gamma
x)(r^2-t^2)_-^{j-\frac{n-3}{2}-2}}{4^j\Gamma(j-\frac{n-3}{2}-1)}dt +
\, O(1)\, \exp(O(T))
\end{multline}
In other words, the contribution of the $0$-th term of the
parametrix  \eqref{param} on the diagonal into the formula
\eqref{klamt:param} cancels out the contribution of the leading
term in  Weyl's  law on the left-hand side of \eqref{pretrace2}.
\end{lemma}
\noindent{\bf Proof.}
Substituting (\ref{disk}) into (\ref{Weyl}) we get:
\begin{equation*}
dN_x(r)=\frac{r^{n-1}}{2^{n-1}\pi^{n/2}\Gamma(n/2)}+dR_x(r),
\end{equation*}
and hence the contribution of the leading term in Weyl's  law
into the left-hand side of (\ref{pretrace2}) is:
\begin{equation}
\label{leadcont}
\int_{-\infty}^\infty \int_{0}^{\infty}
\frac{\psi(t/T)r^{n-1} \cos(r\, t) dr \, dt}{2^{n-1}\pi^{n/2}\Gamma(n/2) T}
\end{equation}
The contribution of the $0$-th term
of the parametrix on the diagonal to (\ref{klamt:param}) is
\begin{equation}
\label{paramcont} \int_{-\infty}^\infty\frac{\psi(t/T) \cos(\lambda
t) |t|^{-n}dt}{\pi^{\frac{n-1}{2}}\Gamma(\frac{1-n}{2})T},
\end{equation}
where $\frac{|t|^{-n}}{\Gamma(\frac{1-n}{2})}$ is  understood as a
generalized function in the sense of (\cite{Gelfand}). If $n=2m$
is even, $|t|^{-2m}=t^{-2m}$ in the sense of generalized functions
and (see \cite{Zelditch}, \cite{Gelfand}):
\begin{equation}
\label{nchet} t^{-2m}=\operatorname{Re} (t+i0)^{-2m} =
\frac{(-1)^m}{(2m-1)!}  \int_0^{\infty} r^{2m-1}\cos(rt) dr.
\end{equation}
Substituting this into (\ref{paramcont}) and taking into account that
$$\Gamma(\frac{1}{2}-m)=\frac{(-1)^m \, m! \, 4^m \sqrt{\pi}}{(2m)!}$$
one obtains equality between (\ref{paramcont}) and (\ref{leadcont}).
If $n=2m+1$ is odd, we have (\cite{Gelfand}):
\begin{equation}
\label{nnechet}
\frac{|t|^{-n}}{\Gamma(\frac{1-n}{2})}=\frac{(-1)^m \,
\delta^{(2m)}(t)\, m!}{(2m)!}= \frac{m!}{\pi \,
(2m)!}\int_0^{\infty} r^{2m} \cos (rt) dr,
\end{equation}
where the second equality follows by taking the inverse  Fourier
transform of $\delta^{(2m)}(t)$. Substituting (\ref{nnechet}) into
(\ref{paramcont}) we again get equality between (\ref{paramcont})
and (\ref{leadcont}). Formula (\ref{pretrem}) is then a
consequence of (\ref{klamt:param}) and (\ref{pretrace2}). This
completes the proof of the lemma. \qed

\smallskip

Lemma \ref{diag:claim} implies that
\begin{equation}\label{pretraceon}
\tilde{k}_{\lambda,T}(x)=\sum_{\omega\in\Gamma\setminus{\rm Id}}
K_{\lambda,T}(x,\gamma x)+ \widetilde{K}_{\lambda,T}(x),
\end{equation}
where
\begin{equation}
\label{lll}
\widetilde{K}_{\lambda,T}(x)=K_{\lambda,T}(x,x)-\int_{-\infty}^\infty
\frac{\psi(t/T)\cos(\lambda t)|t|^{-n}}{\pi^{\frac{n-1}{2}}
\Gamma(\frac{1-n}{2})T} dt.
\end{equation}
Formula (\ref{pretraceon}) indicates that in order to estimate the
remainder $R_x(\lambda)$ one has to study the leading terms in the
right hand side of (\ref{pretraceon}) as $\lambda \to \infty$. These
are the first non-zero on-diagonal term of the parametrix and the
$0$-th off-diagonal terms of the parametrix. As we shall prove
below, in dimensions $n\le 3$ the $0$-th off-diagonal term provides
the principal contribution. In dimension $n\ge 4$ the lower bound
\eqref{bound:weakon} cannot be improved by this method.


\subsection{Off-diagonal leading term}
We shall  consider the cases of $n$ even and $n$ odd separately. If
the dimension $n=2m$ be even, the formula \eqref{param} becomes
\begin{equation}\label{param:even}
E(t,x,y)=\frac{1}{\pi^{m-1/2}} |t| \sum_{j=0}^{\infty}(-1)^j
u_j(x,y) \frac{(r^2-t^2)_-^{j-m-1/2}} {4^j \Gamma(j-m+1/2)} \mod
C^{\infty}.
\end{equation}
If the dimension $n=2m+1$ is odd,  the formula \eqref{param} reads
\begin{equation}
\label{param:odd1} E(t,x,y)=\frac{|t|}{\pi^m}
\sum_{j=0}^{\infty}(-1)^j u_j(x,y) \frac{(r^2-t^2)_-^{j-m-1}} {4^j
\Gamma(j-m)} \mod C^{\infty}.
\end{equation}
Consider first the off-diagonal leading term in even dimensions.
\begin{lemma}\label{chetoffdiag}
Let $n=2m$. The contribution of the leading ($0$-th)  off-diagonal
term of the parametrix to \eqref{klamt:param} as $\lambda \to
\infty$ is
\begin{equation*}
\int_{r}^\infty \frac{\psi(t/T) \, u_0(x,y) \, t\, \cos(\lambda t)
dt} {\pi^{m-\frac{1}{2}}T\Gamma(\frac{1}{2}-m)\,
(t^2-r^2)^{m+\frac{1}{2}}}=Q_1\frac{\psi(r/T)\lambda^{m-\frac{1}{2}}
\sin(\lambda
r+ \phi_m)}{\sqrt{g(x,y)\, r^{2m-1}}\, T} +
O(\lambda^{m-\frac{3}{2}}),
\end{equation*}
where $Q_1$ is some nonzero constant, $r=d(x,y)$, and $\phi_m=\pi/4$
if $m$ is odd, $\phi_m=3\pi/4$ if $m$ is even.
\end{lemma}
\noindent {\bf Proof.} The orders of the leading term and the
remainder follow from \cite[pp. 267-268]{Berard}. To compute the
leading term explicitly, we rewrite the left-hand side as
\begin{multline*}
\int_{r}^\infty \frac{\psi(t/T) \, u_0(x,y) \, t\, \cos (\lambda t)
dt} {\pi^{m-\frac{1}{2}}T\Gamma(\frac{1}{2}-m)\,
(t^2-r^2)^{m+\frac{1}{2}}}=\\
Q_1' \int_{r}^\infty \frac{\psi(t/T) \,t\, \cos (\lambda t)} {T
\sqrt{g(x,y)} (t+r)^{m+\frac{1}{2}}} d\left(
\frac{1}{(t-r)^{m-\frac{1}{2}}}\right),
\end{multline*}
where $Q_1'$ is some non-zero constant. Let us integrate by parts
$m$ times. In order to get the highest power of $\lambda$, all the
differentiations should be applied to $\cos \lambda t$. Otherwise,
the power of $\lambda$ decreases at least by $1$ which implies the
error estimate (the same argument is used in Lemma
\ref{nechetoffdiag}). Differentiating a cosine $m$ times one gets
either a cosine or a sine, depending on the parity of $m$. This
explains different phases $\phi_m$ for $m$ even and odd. The last
step is to make a change of variables $s=t-r$ and to obtain an
Erdelyi-type integral whose asymptotics is well-known (cf.
\cite[Proposition A.1]{Don}). \qed

We do not compute the constant $Q_1$ (as well as the constants
$Q_2,\dots, Q_6$ defined below) since their explicit value is not
important. It only matters that all these constants are non-zero
which one can easily check.
\begin{remark}
In \cite[Lemma III.2]{K} an equivalent statement is obtained for
$m=1$.
\end{remark}

In order study the contribution of the leading off-diagonal  term of
\eqref{param:odd1} in \eqref{klamt:param} we transform it as
follows:
$$
\begin{aligned}
\frac{(-1)^m u_0(x,y)|t| (r^2-t^2)_-^{-m-1}}{\pi^m\Gamma(-m)}
&=\frac{u_0(x,y)|t|}{\pi^m}(-1)^m\delta^{(m)}(r^2-t^2)\\
=\frac{u_0(x,y)|t|(-1)^m}{\pi^m(r+t)^{m+1}}\delta^{(m)}(r-t) &=
\frac{u_0(x,y)|t|}{\pi^m(r+t)^{m+1}} \delta^{(m)}(t-r)
\end{aligned}
$$
\begin{lemma}
\label{nechetoffdiag} Let $n=2m+1$. The contribution of the leading
($0$-th) off-diagonal term of \eqref{param:odd1} in
\eqref{klamt:param} is equal to
\begin{equation*}
\frac{-u_0(x,y)}{T\pi^m} \int_{-\infty}^\infty
\frac{\psi(t/T)|t|\cos(\lambda t)}{(r+t)^{m+1}} \delta^{(m)}(t-r) dt
=Q_2 \frac{\lambda^{m} \, \psi(r/T)  \sin(\lambda\,r +
\phi_m)}{\sqrt{g(x,y)}\, r^m\, T}+O(\lambda^{m-1}),
\end{equation*}
where $Q_2$ is some nonzero constant, $r=d(x,y)$,  and $\phi_m=0$ if
$m$ is odd, $\phi_m=\pi/2$ if $m$ is even.
\end{lemma}
\noindent{\bf Proof.} This lemma follows from \cite[p. 269]{Berard}.
In order to compute the principal term explicitly  one has to
integrate $m$ times by parts applying all derivations to $\cos
(\lambda t)$. Different phases $\phi_m$ for $m$ even and odd appear
for the same reason as in the proof of Lemma \ref{chetoffdiag}. The
principal term is given by
$$
Q_2' \lambda^m \frac{u_0(x,y)}{T\pi^m} \int_{-\infty}^\infty
\frac{\psi(t/T)|t|\sin(\lambda t+\phi_m)}{(r+t)^{m+1}} \delta(t-r)
dt =Q_2 \frac{\lambda^{m} \, \psi(r/T)  \sin(\lambda\,r +
\phi_m)}{\sqrt{g(x,y)}r^m\, T},
$$
where $Q_2, Q_2'$ are some non-zero constants.
 \qed

The lemmas above imply the following
\begin{prop}\label{leading}
The integral $K_{\lambda,T}(x,y)$ defined by  \eqref{klamt1}
satisfies for any $x\neq y \in M$ as $\lambda\to\infty$:
\begin{equation}
\label{leadterm} K_{\lambda,T}(x,y)=
\frac{Q_3\lambda^{\frac{n-1}{2}}\psi(r/T)}{T \sqrt{g(x,y)\,
r^{n-1}}}\;\sin(\lambda r+\phi_n)\;+\; O(\lambda^{\frac{n-3}{2}}).
\end{equation}
Here $r=d(x,y)$, $\phi_n=\frac{\pi}{4}(3-(n\, \operatorname{mod}
4))$, $Q_3$ is a non-zero constant.
\end{prop}
\noindent {\bf Proof.} The proposition follows from
\eqref{klamt:param} and Lemmas \ref{chetoffdiag},
\ref{nechetoffdiag}. The terms in \eqref{klamt:param} for $j>0$
contribute only to the remainder in \eqref{leadterm}. \qed

\subsection{On-diagonal leading term} The on-diagonal leading term is
provided by the coefficient $u_{\kappa_x}(x,x)=a_{\kappa_x}(x)$,
where $\kappa_x$ is defined as in Theorem \ref{weakon}. If the
scalar curvature $\tau(x) \neq 0$ then $\kappa_x=1$. As for the
off-diagonal term, we consider even and odd dimensions separately.
\begin{lemma}
\label{chetdiag} Let $n=2m$ and  $m\ge \kappa_x+1$. The contribution
of the leading ($j=\kappa_x$) on-diagonal term of the parametrix as
$\lambda \to \infty$ in \eqref{pretraceon} is given by
\begin{equation*}
\int_{-\infty}^\infty \frac{\psi(t/T)\, (-1)^{\kappa_x}\,
t^{2\kappa_x-2m}a_{\kappa_x}(x) \cos(\lambda t) dt}
{4^{\kappa_x}\pi^{m-\frac{1}{2}}\Gamma(\kappa_x-m+\frac{1}{2})\, T}=
Q_4\frac{\lambda^{2m-2\kappa_x-1}}{T}+O(\lambda^{2m-2\kappa_x-3}),
\end{equation*}
where $Q_4$ is some non-zero constant.
\end{lemma}

\noindent{\bf Proof.} This lemma can be extracted from \cite[p. 266,
formula (56)]{Berard}. One should expand $\psi(t/T)$ into Taylor
series near $0$ up to the order $2m-2\kappa_x$. Note that all terms
containing odd powers of $t$ vanish after the integration. Using the
fact that the Fourier transform of the generalized function
$t^{-2\alpha}$ is of order $\lambda^{2\alpha-1}$ (\cite{Gelfand}),
we obtain the principal term and the error estimate. \qed

\smallskip

\begin{lemma}
\label{nechetdiag} Let $n=2m+1$ and $m \ge \kappa_x+1$. The
contribution of the the first ($j=\kappa_x$) on-diagonal term of
\eqref{param:odd1} to \eqref{pretraceon} is
\begin{multline}
\label{ondiag1} \frac{(-1)^m a_{\kappa_x}(x)
(m-\kappa_x)!}{4^{\kappa_x} \pi^m (2m-2\kappa_x)! \,T}
\int_{-\infty}^{\infty}\psi(t/T)\, \cos (\lambda t)\,
\delta^{(2m-2\kappa_x)}(t) dt =\\ Q_5
\frac{\lambda^{2m-2\kappa_x}}{T}+O(\lambda^{2m-2\kappa_x-2}),
\end{multline}
where $Q_5$ is some nonzero constant.
\end{lemma}
\noindent{\bf Proof.} This lemma can be extracted from \cite[p.
268]{Berard}. To get the left hand side of \eqref{ondiag1} from
\eqref{param:odd1} we use the first equality in \eqref{nnechet}. The
asymptotics then  follows from the fact that in order to get the
highest power of $\lambda$,  all $2m-2\kappa_x$ derivations coming
from $\delta^{(2m-2\kappa_x)}(t)$ should be applied to $\cos
(\lambda t)$. Note that only an even number of derivations applied
to $\cos (\lambda t)$ produces a non-zero contribution, hence the
orders of the principal and the error terms differ  by $2$. \qed

As consequence of these two lemmas we have the following
\begin{prop}\label{leadingon}
Assume that $n>2 \kappa_x+1$.  Then the integral
$\widetilde{K}_{\lambda,T}(x)$ defined by \eqref{lll} satisfies as
$\lambda\to\infty$:
\begin{equation}
\label{leadtermon} \widetilde{K}_{\lambda,T}(x)=
\frac{Q_6\lambda^{n-2\kappa_x-1}}{T}\;+\;
O(\lambda^{n-2\kappa_x-3}),
\end{equation}
where $Q_6$ is a non-zero constant.
\end{prop}

\smallskip

\noindent{\bf Proof of Proposition \ref{leadingon}.} The
proposition follows from \eqref{klamt:param} and Lemmas
\ref{chetdiag}, \ref{nechetdiag}. The terms in \eqref{klamt:param}
for $j>\kappa_x$ contribute only to the remainder in
\eqref{leadtermon}. \qed


\subsection{Proof of Theorem \ref{weak}}
\label{35} Let $x,y \in X$ be two arbitrary points, $r$ be the
distance between them. Assume that $x$ and $y$ are not conjugate
along any shortest geodesics joining them (and since the geodesics
are shortest, they contain no conjugate points between $x$ and $y$
as well). In particular, this implies that there is only a finite
number of shortest geodesics joining $x$ and $y$, and,  moreover,
there exists $\epsilon>0$ such that any other geodesic joining $x$
and $y$ has length greater than $r+\epsilon$ \cite{Milnor1}. Take
$T=r+\epsilon$ and consider the formula (\ref{klamt}). The
function $\psi$ is supported on $[-1,1]$, hence  one can
approximate the wave kernel $e(t,x,y)$ in \eqref{klamt}  by
summing up the parametrices \eqref{param} (cf. \cite[Theorem
4.3]{Kannai}) along geodesics of length $r$ only. Since $x$ and
$y$ are not conjugate along such geodesics, the parametrices are
well defined and the difference between $e(t,x,y)$ and the sum of
the parametrices is a smooth function in $t$ for $t\in[0,T]$.
Therefore, we can apply the results of the previous section
working directly on  $X$ (and not on $M$). In particular Lemma
\ref{leading} gives us
\begin{equation}
\label{klamtweak}
 k_{\lambda, T}(x,y)
=Q_7\lambda^{\frac{n-1}{2}}\sin{(\lambda r + \phi_n}) +
O(\lambda^{\frac{n-3}{2}}),
\end{equation}
where $Q_3$ is some nonzero constant. Assume for contradiction
$N_{x,y}(\lambda)=o(\lambda^{\frac{n-1}{2}})$.
Set $a=\frac{n-1}{2}$ in Lemma \ref{above:weak1}.
Compare Lemma \ref{above:weak1}
with \eqref{klamtweak} and fix a small number $\nu>0$.
We obtain a contradiction choosing  a sequence
$\lambda_k\to \infty$,  such that
$|\sin(\lambda_k r + \phi_n)|>\nu$. \qed

\subsection{Proof of the Theorem \ref{weakon} (ii)}
\label{36} The proof is similar to the proof of Theorem
\ref{weak}. Consider the orbit $\Gamma x$ of $x$ in the universal
cover $M$. Let
$$
r=\inf_{\operatorname{Id}\neq\omega\in\Gamma} d(x,\omega x).
$$
Since $x$ is not conjugate to itself along any shortest geodesic
loop, the infimum above is attained for finitely many $\omega$-s,
and there exists $\epsilon >0$ such that there are no points in
$\Gamma x\setminus\{x\}$ whose distance to $x$ lies in
$(r,r+\epsilon]$ (\cite{Milnor1}). Choose $T=r+\epsilon/2$ in
\eqref{klamt1}. Only the closest to $x$ lattice points contribute
to the right hand side of \eqref{pretrace2}. This contribution is
of order $\lambda^{\frac{n-1}{2}}$. By Proposition \ref{leadingon}
the contribution of $\omega=\Id$ is of order
$\lambda^{n-2\kappa_x-1}$. If $n-4\kappa_x+1 <0$,
$\lambda^\frac{n-1}{2}$ dominates and hence by the same argument
as in the proof of Theorem \ref{weak} it provides the lower bound
$R_x(\lambda)=\Omega(\lambda^{\frac{n-1}{2}})$ (moreover, in this
case $\frac{n-1}{2} > n-2\kappa_x$ and hence \eqref{bound:weakon}
yields a weaker bound). This completes the proof of
\eqref{weakncp}. \qed

\smallskip

Clearly, applicability of  \eqref{bound:weakon} and \eqref{weakncp} depends
on the geometry of a manifold. For instance, on flat tori
\eqref{bound:weakon} does not give any information since
$u_j(x,x)=~0$ for all $j\ge 1$. However, tori have no conjugate
points and one can use \eqref{weakncp}.

\subsection{Proof of Theorem \ref{thm:osc}}
\label{pf:osc}
The proof of this result is analogous to the proof of
Theorem \ref{weakon} (ii). The condition $n-4\kappa_x+1 <0$ can be
omitted since the sum in \eqref{osc} cancels out the contributions
of the on-diagonal terms of the parametrix up to the order
$[\frac{n-1}{2}]$ to the pretrace formula. This can be checked by
inspection of the coefficients similarly to the proof of Lemma
\ref{diag:claim}. Contributions of the higher-order on-diagonal
terms are negligible. Therefore, the leading contribution comes from
the $0$-th off-diagonal term which is of order
$\lambda^{\frac{n-1}{2}}$, and it yields the lower bound
$R_x^{osc}(\lambda)=\Omega(\lambda^{\frac{n-1}{2}})$. \qed

\begin{remark} It would be interesting to understand whether one
can omit the non-conjugacy condition in Theorems \ref{weak},
\ref{weakon} and \ref{thm:osc}. The proofs will not work since the
parametrix is not well defined for conjugate points. However, our
lower bounds may still hold. Moreover, estimates may get even
stronger: for instance, at any point of a round sphere
$R_x(\lambda)=\Omega(\lambda^{n-1})$.

\end{remark}

\section{Spectral estimates from below on negatively curved manifolds}\label{dyn}

\subsection{Sums over geodesic segments}
\label{seg}
Let $X$ be a negatively curved manifold with the sectional
curvature satisfying \eqref{curvature},  $M$ be the universal
cover of $X$, and let $x,y\in M$. Since $X$ has no conjugate
points, the parametrix is well-defined for all times on $M$.

Our strategy to strengthen Theorems \ref{weak} and \ref{weakon} is
as follows. Let $T$ grow with $\lambda$ and assume that we can
force all the $\sin(\lambda r + \phi_n)$ in Proposition \ref{leading} be of
the same sign and bounded away from zero (see section
\ref{section:Dirichlet}). Consider the sum over geodesic segments
starting at the point $x$:
\begin{equation}
\label{dynamicsum} S_{x,y}(T)=\sum_{r_\omega \le T,\, \omega\in
\Gamma}\frac{1}{\sqrt{g(x,\omega y)\,r_\omega^{n-1}}},
\end{equation}
where $x \neq \omega y$ for all $\omega \in \Gamma$, $r_\omega=d(x,\omega y)$
and $g=\sqrt{\det
  g_{ij} }$. The idea is to estimate \eqref{dynamicsum} from
below by a function going to infinity as $\lambda \to \infty$.
\begin{theorem}\label{sumgrow}
There exists a constant $C_0>0$ such that
\begin{equation}\label{grow1}
S_{x,y}(T)\ \geq\ \frac{C_0}{T} e^{P\left(-\frac{\cH}{2}\right)\cdot T}
\end{equation}
as $T\to\infty$, where $P$ is the topological pressure
\eqref{pressure} and $\cH$ is the Sinai-Ruelle-Bowen potential
\eqref{SRB}.
\end{theorem}
It was shown in \eqref{press:est} that
$P\left(-\frac{\cH}{2}\right)\geq (n-1)K_2/2$, hence $S_{x,y}(T)$
grows exponentially in T.

In the Appendix we shall prove a stronger version of \eqref{grow1}:
\begin{equation}
\label{grow11}
S_{x,y}(T)\ \geq\ C_0 e^{P\left(-\frac{\cH}{2}\right)\cdot T}
\end{equation}
However, a weaker estimate \eqref{grow1} is sufficient for the proofs of
Theorems \ref{main} and
\ref{main:offdiag}.

The proof of Theorem \ref{sumgrow}
is divided into several steps described below.
\subsection{Jacobi fields and the geodesic flow}
Fix a point $x\in M$ and let $\gamma(t)$ be a minimizing geodesic
from $x$ to $y \in M$ such that $\gamma(0)=x$, $\gamma(r)=y$,
$\gamma'(0)=v$, where $v\in T_xM$ is a unit vector. Choose an
orthonormal basis in $T_xM$ consisting of $v$ and vectors $e_1,
e_2,\dots, e_{n-1}$. To every vector $e_j$, $j=1,\dots, n-1$, there
corresponds a perpendicular Jacobi field $Y(t)$ along the geodesic
$\gamma$ with the initial conditions $Y_j(0)=0, Y'_j(0)=e_j$. Then
(see \cite[p. 27]{Berard})
\begin{equation}
\label{Jacdet} g(x,y)=\frac{\sqrt{\det
  (Y_j(r),Y_k(r)) }}{r^{n-1}}
\end{equation}
Let us recall the construction of the canonical (Sasaki) metric on
the tangent bundle $TM$ (see \cite{E}). Let $\Pi:TM\to M$ be the
projection map and let $\xi=(x,v)\in TM$, $\pi(\xi)=x$. Let $\K:
T_\xi(TM)\to T_xM$ be the connection map. There exists a canonical
splitting $T_\xi TM=T_\xi(T_xM)\oplus \Horiz(\xi)$ into vertical
$\operatorname{Im} \K=T_\xi(T_xM)$ and horizontal
$\Horiz(\xi)=\operatorname{Im}d\pi$ subspaces. For any two vectors $z,w \in
T_\xi(TM)$ the Sasaki metric on $TM$ is defined by
\begin{equation}
\label{Sasaki}
 (z,w)_{T_\xi(TM)}=(d\Pi z, d\Pi w)_{T_xM} + (\K z,
\K w)_{T_xM}.
\end{equation}
The induced metric on the unit tangent bundle $SM$ is also referred to
as Sasaki metric.

Consider the geodesic flow $G^t:TM \to TM$. Then we have
(\cite{E}):
\begin{equation}
\label{Jacflow} Y_j(t)=d\Pi \circ dG^t (\tilde e_j),\,\, Y_j'(t)=\K \circ
dG^t(\tilde e_j),
\end{equation}
where $\tilde e_j=\K^{-1}e_j$, $j=1,\dots,n-1$. From the definition of
the Sasaki metric it follows that since $e_j$ are orthonormal
vectors, $\tilde e_j$ are also orthonormal.

\subsection{Jacobian of the geodesic flow}
Consider now the unit tangent bundle $SM$. On negatively-curved
manifolds there exists another natural splitting of $T_\xi(SM)$ into a
direct sum of $DG^t$-invariant subspaces. It comes from the Anosov property
of the geodesic flow $G^t$ on $SM$:
\begin{equation}\label{split} T_\xi(SM)=E_\xi^u\oplus E_\xi^o\oplus E_\xi^s.
\end{equation} Here $E_\xi^u$ is the {\em unstable} subspace  of dimension
$(n-1)$, $E_\xi^s$ is  the {\em stable} subspace  of dimension $(n-1)$, and
$E_\xi^o$
is a one-dimensional subspace tangent to the flow.
\begin{lemma}\label{lemma:neust}
There exists a universal constant $C_1>0$ such that
\begin{equation}
\label{neust} \sqrt{g(x,y)} r^{\frac{n-1}{2}}<C_1 \det
dG^r|_{E_\xi^u},
\end{equation}
where $\Pi(\xi)=x, r=d(x,y)>0$ and $E_\xi^u$ is the unstable
subspace of $T_\xi(SM)$.
\end{lemma}
\noindent{\bf Proof.} Due to \eqref{Jacdet} and \eqref{Jacflow} it
remains to show that
\begin{equation}
\label{remains}
 \det  \left(d\Pi \circ dG^r(\tilde e_j),d\Pi\circ
dG^r(\tilde e_k)\right) \le C_1^2 (\det dG^r|_{E_\xi^u})^2
\end{equation}
for some constant $C_1>0$. First, we argue that
\begin{equation}
\label{bt} \det  \left(dG^r(\tilde e_j), dG^r(\tilde e_k)\right)_{T_\xi(SM)} \,
\ge \, \det
 \left(d\Pi \circ dG^r(\tilde e_j),d\Pi\circ
dG^r(\tilde e_k)\right) ,
\end{equation}
where the scalar products on the left are taken with respect to
the Sasaki metric \eqref{Sasaki}. Indeed, the matrix on the left
is the sum of the matrix on the right and a symmetric positive
definite matrix (coming from the second term in \eqref{Sasaki})
and hence it has a greater determinant.

Denote $\ver(\xi)=T_\xi(S_xM)$ the vertical component of the canonical
splitting of $SM$.
Note that $ \det\left(dG^r(\tilde e_j), dG^r(\tilde e_k)\right)_{T_\xi(SM)}=
(\Jac_{\ver(\xi)}G^r)^2$, where
$\Jac_{\ver(\xi)} G^r$
is the Jacobian of the map $G^r$ restricted to $\ver(\xi)$ (\cite[p. 607]{HK}).
We need to prove
\begin{equation}
\label{ff}
\frac{\Jac_{\ver(\xi)}G^r}{\Jac_{E^u_\xi}G^r}\le C_1.
\end{equation}
It is sufficient to prove this for an integer $r$. Indeed, a volume
comparison argument implies
$\frac{\Jac_{\ver{\xi}}G^r}{\Jac_{\ver{\xi}}G^{[r]}}\leq c$, where
$c$ is a universal constant depending only on the curvature bounds.
At the same time, it is known that $\Jac_{E_\xi^u}G^t$ is a monotone
increasing function in $t$ and hence $\Jac_{E_\xi^u}G^r \ge
\Jac_{E_\xi^u}G^{[r]}$.
Taking a logarithm in \eqref{ff} and writing  $G^r=G^1\circ\dots\circ G^1$
($r$ times) we get:
\begin{equation*}
|\log \Jac_{\ver(\xi)}G^r - \log \Jac_{E^u_\xi}G^r|\le
\sum_{i=1}^{r-1} |\log \Jac_{dG^i\ver(\xi)}G^1 -
\log \Jac_{E^u_{G^i \xi}}G^1|.
\end{equation*}
Let us show that the distance between the subspaces $dG^i\ver(\xi)$ and
$E^u_{G^i \xi}$
is converging to zero exponentially fast.
Here distance is understood in the following sense
(see \cite[p. 191]{AA}). For any $v \in dG^i\ver(\xi)$ let $v=v_u+v_s$ where
$v_u\in E^u_{G^i \xi}$, $v_s \in E^s_{G^i \xi}$. Set
\begin{equation}
\label{arav}
\operatorname{dist}(dG^i\ver(\xi),E^u_{G^i \xi})=
\sup_{||v||=1,\, v\in dG^i \ver{\xi}} \frac{||v_s||}{||v_u||}
\end{equation}
It then follows from the Anosov property of the geodesic flow and the results
of \cite[p. 456]{E}
that \eqref{arav} is well defined and
$$\operatorname{dist}\left(dG^i\ver(\xi),E^u_{G^i \xi}\right)
\le C_1'e^{-\alpha i}$$
for some positive constants $C_1'$ and $\alpha$.
We now remark that $\log \Jac_V G^1$ depends smoothly on a subspace
$V \subset T_\eta(SM)$,
and therefore
$$
|\log \Jac_{dG^i\ver(\xi)}G^1 -
\log \Jac_{E^u_{G^i \xi}}G^1| \le C_1'' e^{-\alpha i},
$$
where $C_1''$ is some positive constant that can be chosen independently of
$i$ by a compactness
argument ($G^1$ on $SM$ is a lift of $G^1$ from a compact manifold $SX$).
Hence
$$
|\log \Jac_{\ver(\xi)}G^r - \log \Jac_{E^u_\xi}G^r|\le C_1''\sum_{i=0}^{r-1}
e^{-\alpha i}
\le \frac{C_1''}{1-e^{-\alpha}}.
$$
This completes the proof of the lemma.
\qed

\smallskip

As an immediate corollary of Lemma \ref{neust} we get
\begin{equation}
\label{sxyt}
 S_{x,y}(T)\ge \frac{1}{C_1} \sum_{r_\omega \le
T,\, \omega\in \Gamma}\frac{1}{\sqrt{\det
dG^{r_\omega}|_{E_{\xi_\omega}^u}}},
\end{equation}
where $\xi_\omega=(x, v_\omega)\in SM$,
$v_\omega=\gamma'_\omega(0)$, $\gamma_\omega$ is the shortest
geodesic joining $x$ and~$\omega y$.

Denote
\begin{equation}\label{Z:def}
Z(r_\omega,\xi_\omega)=\int_{0}^{r_\omega}\left.
\frac{d}{d\tau}\right|_{\tau=0}
\ln\det\left(dG^{s+\tau}|_{E_{\xi_\omega}^u}\right)\;ds =\ln \det
dG^{r_\omega}|_{E_{\xi_\omega}^u}=\int_0^{r_\omega} \cH(G^s
\xi_\omega)ds,
\end{equation}
where $\cH$ is the Sinai-Ruelle-Bowen potential \eqref{SRB}. Here
we used the invariance of the unstable foliation with respect to
the geodesic flow. Hence,
\begin{equation}
\label{sxytsrb} S_{x,y}(T)\ge \frac{1}{C_1} \sum_{r_\omega \le
T,\, \omega\in \Gamma}
\exp\left(-\frac{1}{2}\int_0^{r_\omega} \cH(G^s \xi_\omega)ds
\right).
\end{equation}

\subsection{From geodesic segments  to closed geodesics}
Our next aim is to estimate from below the sum \eqref{sxyt} taken
over geodesic segments starting at $x$ by a sum over closed
geodesics, and then apply techniques of the thermodynamic
formalism. It is well-known that each conjugacy class $[\omega]
\subset \Gamma$ corresponds to a unique closed geodesic
$\gamma_{[\omega]}$ on $X$ of length $l_{[\omega]}$.
Choose a representative $\omega_c \in [\omega]$ corresponding to
this closed geodesic and assume that
it joins $q \in M$ and $\omega_c q \in M$, such that
$d(x,q) \le D$ and $d(\omega y, \omega_c q) \le D$, where $D$ is the diameter
of $X$.
Set $\xi_{[\omega]}=(q, \gamma_{[\omega]}'(0))\in SM$. To distinguish between
distances
on $M$ and on $SM$ we shall write $d_M$
and $d_{SM}$ respectively. Here $d_{SM}$ is induced by the Sasaki metric.

\begin{lemma}\label{comparable1}
There exists $C_2>0$  such that for any $x,y\in M$ and
$\omega \in\Gamma$, $d(x, \omega y)\geq D$ we have
\begin{equation}\label{compare1}
|Z(r_\omega, \xi_\omega)-Z(l_{[\omega]}, [\xi_{[\omega]}])|\leq
C_2.
\end{equation}
\end{lemma}
Let $\gamma_1=[x_1,y_1]$ and $\gamma_2=[x_2,y_2]$ be two geodesic
segments in $M$ such that $\dist(x_1,x_2)<D$ and $\dist(y_1,y_2)<D$.
We shall prove the following general statement implying
Lemma \ref{comparable1}:
\begin{equation}
\label{comparable2}
|Z(r_1, \xi_1)-Z(r_2,\xi_2)|\leq C_2,
\end{equation}
where $\xi_j=(x_j,\gamma_j'(0))$, $j=1,2$ and $r_j=\dist(x_j,y_j)$.
\begin{remark}
\label{diagonall}
It suffices to prove \eqref{comparable2} for two
geodesics starting at the same point on $M$: apply it first to
$(y_1,z_1)$ and $(y_1,z_2)$, then to $(y_1,z_2)$ and $(y_2,z_2)$.
\end{remark}
Accordingly, let $v_1,v_2\in SM,$ such that $\Pi(v_1)=x_1=x_2=\Pi(v_2)$.
Let $\gamma_{v_1}(t)=\gamma_1(t)$ and $\gamma_{v_2}(t)=\gamma_2(t)$ be
the corresponding geodesics in $M$.

\smallskip

\noindent {\bf Note:} to simplify notations, instead of writing
$(\gamma(s),\gamma'(s))$ for an element of  $SM$ we shall write $\gamma'(s)$.

A key ingredient in the proof
of Lemma \ref{comparable1} is the following
\begin{prop}\label{diverge}
Suppose that ${\dist}(\gamma_1(t),\gamma_2(t))\leq D$. Then there
exist two constants $A_1>0,B_1>0$ such that for $0\leq s\leq t$,
\begin{equation}
\label{prop444}
{\DIST}(\gamma_1'(s),\gamma_2'(s))\leq A_1 e^{B_1(s-t)}
\end{equation}
\end{prop}

\subsection{Proof of Proposition \ref{diverge}}
This proof was
communicated to the authors by D. Dolgopyat.
We first remark
that it follows from \cite[Lemma 2.1]{Schroeder} that there exists
$\alpha_1>0$ such that
\begin{equation}\label{est:down}
{\rm dist}(\gamma_1(s),\gamma_2(s))\leq \alpha_1 e^{K_2(s-t)},
\end{equation}
thus establishing the desired estimate on $M$.

It follows easily from the definition of
${\DIST}$ that
\begin{equation}\label{updown:trivial}
{\DIST}(v_1,v_2)\leq \sqrt{\dist (x_1,x_2)^2+\pi^2}
\end{equation}
By increasing $\alpha_1$ if necessary, we may assume without loss of
generality that \eqref{est:down} holds for $s\leq t+1$.

We remark that
it follows easily from \eqref{est:down} and
\eqref{updown:trivial} that \eqref{prop444} would hold if
$\dist(\gamma_1(s+1),\gamma_2(s+1))\geq 1/3$. Accordingly,
it suffices to establish \eqref{prop444} for such $0<s<t$ that
$\dist(\gamma_1(s+1),\gamma_2(s+1))<1/3$.

Consider now the following four points on $M$:
\begin{equation}
\label{points:1}
P_1=\gamma_1(s), P_2=\gamma_2(s), Q_1=\gamma_1(s+1),
Q_2=\gamma_2(s+1).
\end{equation}
We know that
$$
\begin{aligned}
d_1:=\dist(P_1,P_2)&\leq \alpha_1e^{K_2(s-t)},d_1\leq 1/3;\\
d_2:=\dist(Q_1,Q_2)&\leq \alpha_1e^{K_2(s+1-t)},d_2\leq 1/3;\\
\dist(P_1,Q_1)&=\dist(P_2,Q_2)=1.
\end{aligned}
$$
It follows from the
triangle inequality that
\begin{equation}\label{krokodil1}
1-d_1\leq d_3:=\dist(P_1,Q_2)\leq 1+d_1.
\end{equation}

Denote by $\gamma_3(s),0\leq s\leq d_3$ the geodesic segment
connecting $P_1$ and $Q_2$ on $M$.
Consider now the following vectors in $SM$:
\begin{equation}\label{vectors:1}
w_1=\gamma_1'(s),w_2=\gamma_2'(s),w_3=\gamma_3'(0),
w_4=-\gamma_2'(s+1), w_5=-\gamma_3'(d_3).
\end{equation}

The idea is to estimate $\DIST(w_1,w_2)$ using the triangle
inequality
\begin{equation}\label{triang1}
\DIST(w_1,w_2)\leq\DIST(w_1,w_3)+\DIST(w_2,w_3).
\end{equation}

We first remark that $\DIST(w_1,w_3)$ is equal to the angle at the
vertex $P_1$ in the geodesic triangle $P_1Q_1Q_2$ with sides
$d_1,d_2,d_3$. We compare the triangle $P_1Q_1Q_2$ with the planar triangle
having the same sides.  It follows from $d_2<1/3$, \eqref{krokodil1},
elementary planar trigonometry and comparison theorem that
the angle $\theta=\angle Q_1P_1Q_2$  is less than $\pi/6$.  Hence,
$\sin\theta > 3\theta/\pi$. Therefore we have
\begin{equation}\label{triang2}
\DIST(w_1,w_3)=\theta < \frac{\pi \sin \theta}{3} <
\frac{\pi d_2}{3 d_3} < \frac{\pi d_2}{2} \le \frac{\pi\alpha_1 e^{K_2}}{2}
e^{K_2(s-t)}
\end{equation}
Here we use that $\sin \theta < d_2/d_3$ which follows from
comparison with a planar triangle, and that $d_3>2/3$ which follows from the
triangle
inequality.
It remains to estimate $\DIST(w_2,w_3)=\DIST(-w_2,-w_3)$. Consider
the vector $w_6=-\gamma_3'(d_3-1)$.  We see that
$$
\DIST(-w_2,-w_3)\leq\DIST(-w_2,w_6)+
\DIST(w_6,-w_3)\leq\DIST(-w_2,w_6)+d_1,
$$
so it suffices to estimate $\DIST(-w_2,w_6)$.

Now,
$$
-w_2=G^1(w_4),\ \ w_6=G^1(w_5),
$$
where $G^1$ is the time one map of the geodesic flow $G^t$. The
distance between $w_4$ and $w_5$ is equal to the angle at $Q_2$ of
the geodesic triangle $P_1Q_2P_2$ with sides $d_3,1,d_1$. By
comparing that triangle to the isometric planar triangle, one can
show as in \eqref{triang2} that
\begin{equation}\label{triang3}
\DIST(w_4,w_5)\leq \frac{\pi d_1}{2}\leq \frac{\pi\alpha_1}{2} e^{K_2(s-t)}.
\end{equation}
${\rm inj}(X)$ is

Now, consider the smooth diffeomorphism $G^1$ restricted to the unit
tangent space of the
quadrilateral $P_1Q_1P_2Q_2$. By compactness it increases distances by at
most a factor
$\alpha_2$.  It follows from \eqref{triang3} that
\begin{equation}\label{triang4}
\DIST(-w_2,w_6)\leq \alpha_2 \DIST(w_4,w_5)\leq
\frac{3\pi\alpha_1\alpha_2}{4} e^{K_2(s-t)}.
\end{equation}

The inequality \eqref{prop444} now follows from
\eqref{triang1}, \eqref{triang2} and \eqref{triang4}.
This finishes the proof of the proposition.
\qed

\subsection{Proof of Lemma \ref{comparable1}}
By Remark \ref{diagonall}, in order to prove \eqref{comparable2}
it suffices to consider two vectors $\xi_1,\xi_2\in SM$ with
$\pi(\xi_1)=\pi(\xi_2)$. Let $\gamma_1$ and $\gamma_2$ be geodesics
as in Lemma \ref{diverge}: say, $\gamma_1=[x,\omega y]$
and $\gamma_2=[x,\omega_c q]$.
Let $l_1, l_2$ be the lengths of $\gamma_1,\gamma_2$
respectively. It follows from the triangle inequality that
$|l_1-l_2|\leq D$.
It thus suffices to estimate the difference
$$
|Z(r_1,\xi_1)-Z(r_2,\xi_2)|=\left|\int_0^{r_1}\cH(G^s\xi_1)ds-
\int_0^{r_2}\cH(G^s\xi_2)ds\right|
$$
Since $\cH$ is uniformly bounded on $SX$ by compactness, and thus also on $SM$,
the inequality \eqref{comparable2} will remain true (with a
different constant) if we increase or decrease the length of
$\gamma_j$ by a uniformly bounded amount.  Since $|l_1-l_2|\leq
D$, we may thus assume without loss of generality that
$\gamma_j$-s have the same length $t$.

Accordingly, to prove \eqref{compare1}, it suffices to show that
there exist $B_2>0$ such that for $t\geq D$,
\begin{equation}\label{comparelog}
\left|\int_0^t [\cH(G^s\xi_1))-\cH(G^s\xi_2)]\; ds\right| \leq
B_2.
\end{equation}

It follows from the H\"older continuity of $H$ (\cite{BR},\cite{Sinai2})
that there exist two constants $A_2>0,\beta>0$ such
that for any $v_1,v_2\in SM$,
$$
|\cH(v_1)-\cH(v_2)|\leq A_2\cdot \DIST(v_1,v_2)^\beta.
$$

It follows from Proposition \ref{diverge} that there exist two
constants $A_1, B_1>0$, such that for $0\leq s\leq t$,
$$
\DIST(G^s\xi_1,G^s\xi_2)\leq A_1 e^{B_1(s-t)}
$$

Call the integral in \eqref{comparelog} $I(t)$.  Now,
$$
\begin{aligned}
|I(t)|&\leq \int_0^t|\cH(G^s\xi_1)-\cH(G^s\xi_2)|ds \leq
A_2\int_0^t\;
\DIST(G^s\xi_1,G^s\xi_2)^\beta\;ds\\
\;&\leq A_2A_1^\beta\;\int_{0}^t e^{B_1\beta(s-t)}\; ds \leq
\frac{A_2A_1^\beta}{B_1\beta}.
\end{aligned}
$$
This proves \eqref{comparelog} with $B_2=A_2A_1^\beta/(B_1\beta)$.
The proof of Lemma \ref{comparable1} is complete.
\qed


\subsection{Proof of Theorem \ref{sumgrow}}
\label{therm}
By Lemma \ref{comparable1} and \eqref{sxytsrb} we have
\begin{equation}
\label{below2}
S_{x,y}(T)\geq C_3\sum_{\omega\in\Gamma: \, 3D\leq
l_{[\omega]}\leq T-2D} \exp\left(-\frac{Z(l_{[\omega]},\xi_{[\omega]})}{2}
\right)
\end{equation}
for some constant $C_3>0$. We sum over the interval $3D\leq l_{[\omega]}\leq T-2D$
since by triangle inequality the difference between the length of the segment
and the corresponding
closed geodesic is at most $2D$, and due to the condition of Lemma
\ref{comparable1} we consider
segments of length $\ge D$. It is a relatively rough bound. Indeed,
we take into account the
contribution of just one element from every conjugacy class. In the
Appendix we refine this estimate
taking into account the number of elements in each conjugacy class.

We now apply quite a deep fact on the equidistribution of closed geodesics
for Anosov flows
proved by methods of thermodynamic formalism.
Using \eqref{below2} we write:
\begin{equation}
\label{below3}
S_{x,y}(T)\geq \frac{C_3}{T}\sum_{\omega\in\Gamma: \, 3D\leq
l_{[\omega]}\leq T-2D} l_{[\omega]} \exp\left(-\frac{Z(l_{[\omega]},
\xi_{[\omega]})}{2}
\right)
\end{equation}
On the other hand, the sum \eqref{below3} over closed geodesics satisfies
the following
asymptotic relation for $T\to \infty$: (see \cite{Par},
\cite[(7.1)]{PP}, \cite[p. 109]{MS}; cf. \cite{Ruelle}),
\begin{equation}\label{dynam:asymp}
\sum_{l_{[\omega]}\leq T} l_{[\omega]}
\exp\left(-\frac{Z(l_{[\omega]},\xi_{[\omega]})}{2} \right)=
\frac{C_4}{P(-\cH/2)}\;e^{P(-\cH/2)T}\left(1+o(1)\right),
\end{equation}
where $$C_4=\int_{SX}\;d\mu_{-\cH/2}.$$
Here $\mu_{-\cH/2}$ is the
{\em equilibrium state} for $-\cH/2$ and $P(-\cH/2)$ is the pressure
defined by \eqref{pressure}. As $T\to \infty$, the contribution of
closed geodesics of lengths $l_{[\omega]}\le 3D$ and $T-2D\le
l_{[\omega]}\le T$ to the sum in \eqref{dynam:asymp} is bounded by a
constant. Therefore, \eqref{dynam:asymp} implies \eqref{grow1} and
this completes the proof of Theorem \ref{sumgrow}. \qed


\section{Proof of the main results.}\label{proofs}

\subsection{Proof of Theorem \ref{main:offdiag} for
$n\not\equiv 3({\rm mod}\, 4)$.}\label{off:main}

Given $\delta>0$, let
$$\alpha=\frac{(1-\delta)P(-\cH/2)}{h}.$$
Assume for
contradiction that the bound \eqref{bound:offdiag} doesn't hold
and thus for some $x,y\in X, x\neq y,$ the off-diagonal spectral
counting function $N_{x,y}(\lambda)$ satisfies the following upper
bound:
\begin{equation}\label{contr:offdiag}
N_{x,y}(\lambda)=O\left(\lambda^{(n-1)/2}
(\log\lambda)^\alpha\right).
\end{equation}

By Lemma \ref{above2} we conclude that  there
exists a constant $C_5>0$ such that for any $\lambda$ we have
\begin{equation}\label{off:contr2}
|k_{\lambda, T}(x,y)|<C_5\lambda^{(n-1)/2}(\log\lambda)^\alpha.
\end{equation}

Using \eqref{pretrace2} we write
$$
k_{\lambda, T}(x,y)= \sum_{\omega\in\Gamma}K_{\lambda,T}(x,\omega y).
$$
We want to estimate the right-hand side from below.
It follows from Proposition \ref{leading} that
\begin{equation}
\label{off:sum1}
k_{\lambda, T}(x,y)= \sum_{\omega\in\Gamma, r_\omega \le T}
\frac{Q\lambda^{\frac{n-1}{2}}}{T}\frac{\,\psi(\frac{r_\omega}{T})}
{\sqrt{g(x,\omega y)\,r_\omega^{n-1}}}\;\sin(\lambda r_\omega
+\phi_n)+O(\lambda^{\frac{n-3}{2}})\exp(O(T)).
\end{equation}
The sum in \eqref{off:sum1} can be taken over $r_\omega\leq T$ since
$\supp\psi=[-1,1]$.
By results of Margulis (\cite{MS}), the number of summands in \eqref{off:sum1} is
less than $C_6 e^{hT}$ for some positive constant $C_6$. Since the error
for each term is
$O(\lambda^{(n-3)/2})$ by Proposition \ref{leading}, for the whole sum the error
is $O(\lambda^{\frac{n-3}{2}})\exp(O(T))$.

\subsection{Dirichlet box principle}
\label{section:Dirichlet}
We want all the terms
$\sin(\lambda r_\omega+\phi_n)$ in \eqref{off:sum1} to have the
same sign, and $|\sin(\lambda r_\omega+\phi_n)|$ to be bounded from below
by some positive constant. Let $\{r_1,r_2,\ldots,r_{N_1}\}$ be all
the distinct values of $r_\omega$ appearing in
\eqref{off:sum1}. It suffices to choose $\lambda$ so that for all
$1\leq j\leq N$,
\begin{equation}\label{dioph}
|e^{i\lambda r_j} -1| <1/10.
\end{equation}
In that case all the angles
$\lambda r_j +\phi_n$ are close to
$$
\phi_n=\pi/4(3-n\;{\rm mod} 4).
$$
Since $n\not\equiv 3({\rm mod} 4)$,  $\sin(\lambda r_j +\phi_n)\approx
\sin \phi_n$
all have the same signs and $|\sin(\lambda r_j +\phi_n)| \ge C_7$
for some positive constant $C_7$.

To establish \eqref{dioph} we apply the following Lemma from \cite{PR}
(a similar approach was also used in \cite{RS} to estimate error terms
from below):
\begin{lemma}\label{box}
Given $r_1,r_2, \dots, r_{N_1}$ and
$M_1>0,Y>1$ there exists $\lambda\in[M_1,M_1Y^{N_1}]$ such that for
all $1\leq j\leq N_1$,
$$
|e^{i\lambda r_j} -1|<1/Y.
$$
\end{lemma}

By Lemma \ref{box}, given $M_1>0$ we can choose
$M_1\leq\lambda\leq M_1\cdot 10^{N_1}$ such that \eqref{dioph} will
hold for all $1\leq j\leq N_1$. We recall that $N_1<C_6 e^{hT}$.

Moreover, given $\delta>0$ we can choose $\psi$ so that
$\psi(x)\geq 1/2$ for $|x|\leq 1-\delta/2$.  Accordingly, if the
radius $r_j$ in \eqref{off:sum1} satisfies
$r_j\leq T(1-\delta/2),$
then $\psi(r_j/T) \ge 1/2$.
Therefore, given $M_1$, there exists $\lambda$ such that
\begin{equation}\label{N:bound}
\ln M_1\leq\ln \lambda\leq \ln M_1+(\ln 10)C_6 e^{hT}.
\end{equation}
It then follows from \eqref{off:sum1} that there exists $A_3>0$
such that
\begin{equation}\label{off:sum2}
|k_{\lambda, T}(x,y)|\geq \frac{A_3\lambda^{(n-1)/2}}{T}
\sum_{\omega:\, \dist(x,\omega
y)\leq(1-\delta/2)T}\frac{1}{\sqrt{g(x,\omega y)r_\omega^{n-1}}}+
O(\lambda^{\frac{n-3}{2}})\exp(O(T))
\end{equation}
In the sequel, we shall choose $T=O(\ln\ln\lambda)$.  It follows
that the error term in \eqref{off:sum2} is
\begin{equation}\label{errorest}
O\left(\lambda^{\frac{n-3}{2}}(\log\lambda)^{O(1)}\right), \qquad
{\rm provided}\ T=O(\ln\ln\lambda).
\end{equation}

Using the estimate \eqref{grow1} proved in Theorem \ref{sumgrow}
we conclude from \eqref{off:sum2} and \eqref{errorest} that
\begin{equation}\label{off:sum3}
|k_{\lambda, T}(x,y)|\geq
\frac{B_3}{T^2}\lambda^{\frac{n-1}{2}}e^{P(-\cH/2)(1-\frac{\delta}{2})T},
\end{equation}
for some constant $B_3>0$.

To obtain contradiction with \eqref{off:contr2}, we should find
$\lambda>0,T>0$ satisfying $\lambda\geq M_1$ and \eqref{N:bound}
$$
C_5\lambda^{(n-1)/2}(\log\lambda)^\alpha\leq
\frac{B_3}{T^2}\lambda^{(n-1)/2}e^{P(-\cH/2)(1-\frac{\delta}{2})T}.
$$
This translates to
\begin{equation}\label{ineq1}
\begin{aligned}
\ln\ln\lambda &\leq\frac{1}{\alpha}\left[T\,P(-\cH/2)(1-\delta/2) +
\ln B_3 -\ln C_5 -2\ln T\right]\\
\; &= \frac{h(1-\delta/2)}{1-\delta} T + \frac{1}{\alpha}\left[
\ln B_3 -\ln C_5 -2\ln T\right]
\end{aligned}
\end{equation}

We can rewrite \eqref{N:bound} as
\begin{equation}\label{ineq2}
\ln\ln M_1\leq\ln\ln\lambda\leq hT+\ln\ln M_1+\ln C_6+\ln\ln 10.
\end{equation}
By Lemma \ref{box}, for any $M_1$  we can find some $\lambda$
satisfying \eqref{ineq2} such that the estimate \eqref{off:sum2}
holds.

Let $\beta$ be a small constant satisfying
$$
0<\beta<\frac{\delta/2}{1-\delta}.
$$
We choose $M_1$ in \eqref{ineq2} to so that
$$
\ln\ln M_1=h\beta T.
$$
Then \eqref{ineq2} becomes
\begin{equation}\label{ineq3}
h\beta T\leq\ln\ln\lambda\leq h(1+\beta)T+\ln C_6+\ln\ln 10
\end{equation}
It follows that $T \approx \frac{1}{h} \ln\ln\lambda$, and thus
\eqref{errorest} holds,
implying \eqref{off:sum3}.

To prove Theorem \ref{main:offdiag}, it suffices to show that if
we choose $T$ large enough, then \eqref{ineq3} would {\em imply}
\eqref{ineq1}.  Indeed, in that case {\em any} choice of $\lambda$
satisfying \eqref{ineq3} would automatically contradict
\eqref{contr:offdiag}, and {\em some} such choice exists by Lemma
\ref{box}.

The largest terms in both \eqref{ineq3} and \eqref{ineq1} are
linear in $T$, so it suffices to compare the coefficients of $T$
in those inequalities. The coefficient in \eqref{ineq1} is equal
to $\frac{h(1-\delta/2)}{1-\delta}$ while the coefficient in
\eqref{ineq3} is equal to $h(1+\beta)$.  By the choice of $\beta$,
the coefficient of $T$ in \eqref{ineq3} is {\em smaller} than that
in \eqref{ineq1}. It follows that for large $T$, the right-hand
side of \eqref{ineq1} is {\em larger} than the right-hand side of
\eqref{ineq3}, therefore \eqref{ineq3} implies \eqref{ineq1},
finishing the proof of Theorem \ref{main:offdiag} for $n\not\equiv
3({\rm mod} 4)$.
\qed

\subsection{Proof of Theorem \ref{main:offdiag} for
$n\equiv 3({\rm mod}\, 4)$.}\label{sec:n3}

The only difference between the proofs in cases $n\not\equiv
3({\rm mod}\, 4)$ and $n\equiv 3({\rm mod}\, 4)$ is that in the latter
case, the angle $\phi_n$ in Proposition \ref{leading} is equal to zero.
Accordingly, the leading term in Proposition \ref{leading} is
proportional to $\sin(\lambda r_\omega)$ and thus we have to ensure
that $\lambda r_\omega$ don't get too close to $\pi\zed$.  In addition,
we want to choose $\lambda$ so that $\sin \lambda r_\omega$  all have the
same sign.   Accordingly, we cannot use
Lemma \ref{box} directly and have to replace it by a different
statement.

The idea of the modification is as follows: the biggest contribution
to the sum \eqref{grow1} in Theorem \ref{sumgrow} comes from the
terms $\omega \in\Gamma$ where $r_\omega=\dist(x,\omega y)$ is
large. Indeed, for any $t>0$ every term in the sum $S_{x,y}(t)$ is
not greater than $1$, therefore the number of terms gives a trivial
upper bound: $S_{x,y}(t) \le C_6 e^{ht}$. Take
\begin{equation}\label{A:ineq}
A_4>\frac{h}{P(-\cH/2)(1-\delta/2)}.
\end{equation}
and write $S_{x,y}(T)=S_1(T)+S_2(T)$, where $S_1(T)$ is the sum
taken over all $\omega$ such that
$r_\omega < T/A_4$ and $S_2(T)$ corresponds to
\begin{equation}\label{dist:large}
T/A_4 \leq r_\omega \leq T.
\end{equation}
Comparing the exponents we deduce that $S_1(T)$ grows slower than
$S_{x,y}(T)$. Indeed, $P(-\cH/2)(1-\delta/2)T>hT/A_4$ by
\eqref{A:ineq}. Therefore, for $T$ large the contribution of
$S_1(T)$ to $S_{x,y}(T)$ is negligible in the sense that
$S_{x,y}(T)=S_2(T)(1+o(1))$.

Let $\{r_1,r_2,\ldots,r_{N_2}\}$ be all the $N_2=O(e^{hT})$ distinct
values of $r_\omega$ appearing in $S_2(T)$. Choose $Y$ satisfying
\begin{equation}\label{B:ineq}
Y>\frac{3A_4}{\pi}.
\end{equation}
By Lemma \ref{box}, given $M_2>0$ we can find
$\mu_1\in[M_2,M_2 Y^{N_2}]$ such that
\begin{equation}\label{dist:init}
\operatorname{dist}(\mu_1 r_j, 2\pi \zed)<1/Y, \qquad 1\leq j\leq N_2.
\end{equation}
The reason for choosing $Y>3A_4/\pi$ will be explained below.

It suffices to find $\lambda$ such that $\sin \lambda r_\omega$ all
have the same sign for $T/A_4 \le r_\omega \le T$, and to
estimate all those $\sin\lambda r_\omega$ from below.  Let
$2\pi k_j$ be the closest multiple of $2\pi$ to $\mu_1 r_j$ in
\eqref{dist:init}, and let
$$
b_j=\mu_1 r_j-2\pi k_j.
$$
We shall prove the following
\begin{lemma}\label{sin:large}
Let $T/A_4\leq
r_1<r_2<\ldots<r_{N_2}\leq T$, $T\gg 1$,  and let $|b_j|\le 1/Y$,
$j=1,\dots, N_2$.
Then there exists $\mu\in[0,1]$ such that
\begin{equation}\label{big:sine}
\sin(\mu r_j+b_j)\geq \frac{1}{2T},\qquad 1\leq j\leq N_2.
\end{equation}
\end{lemma}
\noindent{\bf Proof of Lemma \ref{big:sine}.} We want to find
$\mu\in[0,1]$ such that
\begin{equation}\label{sin:ineq1}
\frac{1}{T}\leq \mu r_j+b_j\leq \pi- \frac{1}{T}, \qquad 1\leq j\leq m.
\end{equation}
The equality \eqref{sin:ineq1} is equivalent to
$$
\mu\in I_j:=\left[\frac{1}{Tr_j}-\frac{b_j}{r_j},
\frac{\pi}{r_j}-\frac{1}{Tr_j}-\frac{b_j}{r_j}\right].
$$
Thus, it suffices to show that
\begin{equation}\label{sin:inter1}
\emptyset\neq\cap_{j=1}^m I_j\cap[0,1].
\end{equation}

The maximum $P_1$ of the left endpoints of $I_j$ is given by
$$
P_1=\max_j\left\{\frac{1}{Tr_j}-\frac{b_j}{r_j}\right\} \leq
\frac{A_4}{T^2}+\frac{A_4}{YT}.
$$
On the other hand, the minimum $P_2$ of the right endpoints of
$I_j$ is given by
$$
P_2=\min_j
\left\{\frac{\pi}{r_j}-\frac{1}{Tr_j}-\frac{b_j}{r_j}\right\}
\geq \frac{\pi}{T}-\frac{A_4}{T^2}-\frac{A_4}{Y T}.
$$
To prove \eqref{sin:inter1} it suffices to show that $P_1<P_2$.
Accordingly, it suffices to show that
$$
\frac{A_4}{T^2}+\frac{A_4}{YT}<\frac{\pi}{T}-\frac{A_4}{T^2}-
\frac{A_4}{Y T}.
$$
Now, it follows easily that the above inequality holds if we
choose $Y>3A_4/\pi$ (as in \eqref{B:ineq}) and $T$ large enough.
It follows that \eqref{sin:inter1} holds, finishing the proof of
Lemma \ref{sin:large}.
\qed

The rest of the proof of the Theorem \ref{main:offdiag} for $n\equiv
3 (\operatorname{mod} 4)$ is the same as in the case $n\not\equiv
3({\rm mod}\, 4)$. The only difference is that in the analogue of
\eqref{off:sum3} we get an extra $T$ in the denominator due to the
dependence on $T$ in \eqref{big:sine}. However, this does not affect
the final result since $T=O(\log \log \lambda)$ and $1/T^3$ is
absorbed by the factor  $-\delta$ in the  power of $\log \lambda$ in
\eqref{bound:offdiag}. \qed

\subsection{Proof of Theorem \ref{main}}\label{on:small}
 The rate of growth of the
error term on the diagonal is determined by the exponents of the
leading terms in Propositions \ref{leading} (for $\omega \neq
\operatorname{Id}$) and \ref{leadingon} (for
$\omega=\operatorname{Id}$).   For $n=2$ and $n=3$ the leading terms
in Proposition \ref{leading} are of order $\lambda^{(n-1)/2}$  and
grow faster than then $\widetilde{K}_{\lambda, T}(x)$.

Note that the geodesic segments considered in section \ref{seg}
in the case of Theorem \ref{main} are geodesic loops. Consequently, the sum
\eqref{dynamicsum} should be taken over $\omega \in
\Gamma\setminus\operatorname{Id}$. We remark that the geodesic {\it loops},
and not the closed geodesics (though they are used in the proof),
contribute to \eqref{dynamicsum} and hence to the growth of
$R_x(\lambda)$ (cf. section \ref{discussion}, \cite[p.175]{Bogom}).

The rest of the proof of Theorem \ref{main} is identical to that of
Theorem \ref{main:offdiag}. In the case $n=3$ the application of the
Dirichlet principle is used appropriately as in section
\ref{sec:n3}.

For $n=2$ and $n=3$ \eqref{bound} gives a better lower bound than
\eqref{bound:weakon}. Note that for $n=3$ the gain is due to the
logarithmic improvement guaranteed by Theorem \ref{sumgrow}. For
$n\ge 4$ the ``apriori'' bound $R_x(\lambda)=\Omega(\lambda^{n-2})$
of \eqref{bound:weakon} can not be improved. This completes the
proof of \eqref{bound}.

It remains to prove \eqref{oscneg}.  We note that just as in the proof of
Theorem \ref{thm:osc}, the sum in \eqref{osc} cancels out the contributions
of the on-diagonal ($\omega=Id$) terms of the wave parametrix up to the order
$[\frac{n-1}{2}]$ to the pretrace formula.
The estimate \eqref{oscneg} now follows from the first formula in
\eqref{bound} exactly in the same way as Theorem \ref{thm:osc} \, is
deduced from Theorem \ref{weakon}, see section \ref{pf:osc}. \qed

\smallskip

\noindent{\bf Acknowledgements.} The authors are  grateful to D.
Dolgopyat for explaining them hyperbolic dynamics leading to the
results in section 4. We would like to thank Yu. Safarov for
pointing out that the bound \eqref{bound:weakon} in the earlier
version of the paper could be improved. We would also like to
thank V. Ivrii, M. Jakobson, V. Jaksic, Y. Kannai, Ya. Pesin, M.
Pollicott, L. Polterovich, P. Sarnak, R. Sharp, A. Shnirelman, M.
Shubin, J. Toth and D. Wise for useful discussions. This paper was
completed while the first author visited IHES and Max Planck Institute
for Mathematics, and their hospitality is greatly appreciated.

\section{Appendix: proof of \eqref{grow11}}\label{betterdyn}
\subsection{Refining the estimate}
\label{group}
In this section we shall prove the following
\begin{prop}\label{prop:grow2}
Let $S_{x,y}(T)$ be as in \eqref{dynamicsum}.  Then
$$
S_{x,y}(T)\ \geq\ C_0 e^{P\left(-\frac{\cH}{2}\right)\cdot T}
$$
\end{prop}

To prove Proposition \ref{prop:grow2}, we group the terms in
\eqref{below2} according to conjugacy classes of $\omega \in\Gamma$.  It
is well-known that each conjugacy class corresponds to a unique
closed geodesic on $X$. This correspondence can be defined as
follows. The universal cover $M$ of $X$ can be tiled by identical
copies of a fundamental domain $\Theta$ for $X$. Those copies are
indexed by elements of
$\Gamma=\pi_1(X)$. Indeed, pick a copy of $\Theta\in S$ and call it
$\Theta_e$, where $e$ is the identity element in $\Gamma$. Then for
any other copy $\widetilde{\Theta}$ there exists a unique
$\omega\in\Gamma$ such that $\omega\Theta_e=\widetilde{\Theta};$ we
then let $\widetilde{\Theta}=\Theta_\omega$.  We thus have

$$
\omega h^{-1}\Theta_h=\Theta_\omega, \qquad\forall \omega,
h\in\Gamma.
$$

Given a {\em simple} closed geodesic $\gamma$ on $X$, choose its
lift to $M$ passing through $\Theta_e$ (which we shall also call
$\gamma$).  Choose a basepoint $z\in\gamma\cap\Theta_e$, and let
$\gamma=[z,\omega z],\omega\in\Gamma$;  this defines $\omega$.
Now, $\gamma$ intersects other copies of the fundamental domain
besides $\Theta_e$; denote those copies by
\begin{equation}\label{intersect1}
\Theta_{h_1},\Theta_{h_2},\ldots,\Theta_{h_n}
\end{equation}
Here $h_1=e,h_n=\omega$ and $\Theta_{h_i}$-s are numbered
consecutively along $\gamma$.  One can show \cite{Milnor} that
$\Gamma$ is generated by
$\{h:\partial\Theta_h\cap\partial\Theta_e\neq\emptyset\}$. Let
$\{a_1,a_1^{-1},\ldots,a_k,a_k^{-1}\}$ denote a set of generators.
Then each of the group elements $h_{i+1}h_i^{-1},1\leq i\leq n-1$
is one of the generators (it maps $\Theta_{h_i}$ to
$\Theta_{h_{i+1}}$).  We set
$$
a_{j_i}:=h_{i+1}h_i^{-1}.
$$
Thus to a geodesic $\gamma$ we have associated an element $\omega$
and a word $w(\gamma)$ in $a_j$-s given by
\begin{equation}\label{word1}
a_{j_n}a_{j_{n-1}}\ldots a_{j_2}a_{j_1}=\omega
\end{equation}

It is clear that
\begin{equation}\label{wordlength1}
n\geq l(\gamma)/{\rm diam}(\Theta),
\end{equation}
where $l(\gamma)$ is the length of the closed geodesic $\gamma$. In
fact, it is shown in \cite[Lemma 2]{Milnor} that (for a fixed set of
generators) the ratio $l(\gamma)/L_\Gamma(w(\gamma))$, the
denominator being the word length of $w(\gamma)$,  is uniformly
bounded above and below. However,  we need only one side of this
estimate, namely inequality \eqref{wordlength1}.

If we choose a basepoint $z$ on $\gamma$ not in $\Theta_e$ but in
adjacent fundamental domain $\Theta_{h_1}$, then the corresponding
$\omega\in\Gamma$ becomes $\tilde{\omega}=h_1\omega h_1^{-1}$, and
the word corresponding to $\tilde{\omega}$ is equal to (cf.
\eqref{word1})
$$
a_{j_1}a_{j_n}a_{j_{n-1}}\ldots a_{j_2}
=a_{j_1}w(\gamma)a_{j_1}^{-1}
$$
is a {\em cyclic shift} of $w(\gamma)$.
Among all such cyclic
shifts, choose one (call it $w_1(\gamma)$) which will have the {\em
smallest} word length (such shift is called a {\em cyclically
reduced} form of $w(\gamma)$). We shall prove the following
\begin{prop}\label{distinct}
All cyclic shifts of a {\em primitive cyclically reduced}
$w_1=w_1(\gamma)$ are different elements of $\Gamma$.
\end{prop}

{\bf Proof of Proposition \ref{distinct}}.  It was shown by
Preissman in \cite{Preis} that any nontrivial commutative subgroup
of $\Gamma$ is infinite cyclic. Let us now assume for contradiction
that two cyclic shifts of $w_1$ coincide. We can assume without loss
of generality that one shift is trivial and equals $w_1$ itself.
So suppose that
$$
w_1=u_1u_2=u_2u_1
$$
for two nontrivial words $u_1,u_2\in\Gamma$.  Then
$u_1$ and $u_2$ generate a nontrivial commutative subgroup of $\Gamma$.
Therefore, by Preissman's theorem it is cyclic and there exists
$u_3\in\Gamma$ such
that $u_1=u_3^k,u_2=u_3^l$ for some $k,l\in\zed$.  We claim that $k$
and $l$ have the same sign, otherwise the word $w_1=u_1u_2$ has
cancellations and is not cyclically reduced.  It follows that
$w_1=u_3^{k+l}$ and is therefore not primitive, which is a
contradiction.  This finishes the proof of the Proposition
\ref{distinct}. \qed

\subsection{Proof of Proposition \ref{prop:grow2}}
We would like to get a better estimate for \eqref{dynamicsum} than
the one proved  in section \ref{therm}.  To do that, consider a
conjugacy class in $\Gamma$ and the corresponding closed geodesic
$\gamma$.  We associate to $\gamma$ a word $w(\gamma)$ of length
$L_\Gamma(w(\gamma))$ as in section \ref{group}.  By Proposition
\ref{distinct}, to each primitive $\gamma$, there  correspond
$l_\Gamma(w(\gamma))$ cyclic shifts of this word, and by
\eqref{wordlength1}
$$
L_\Gamma (w(\gamma))\geq\frac{l(\gamma)}{{\rm diam}(\Theta)}.
$$
It is now clear from the argument of section \ref{therm}, that to
each primitive closed geodesic $\gamma$ we can associate at least
$L_\Gamma(w(\gamma))$ segments connecting $x$ to a point on the
orbit $\Gamma y$. The conclusion of Lemma \ref{comparable1} applies
to each of these segments. Accordingly, the estimate \eqref{below3}
can be improved to
$$
S_{x,y}(T)\geq \frac{C_3}{{\rm diam}(\Theta)}\sum_{\omega\in\Gamma:
\, 3D\leq l_{[\omega]}\leq T-2D} l_{[\omega]}
\exp\left(-\frac{Z(l_{[\omega]},\xi_{[\omega]})}{2} \right),
$$
where $l_{[\omega]}=l(\gamma)$ is the length of the closed geodesic
$\gamma$ corresponding to the group element $\omega$.

An application of \eqref{dynam:asymp} would complete the proof of
Proposition \ref{prop:grow2}, provided one can show that the sum
over {\em primitive} closed geodesics is comparable to the ``total''
sum \eqref{dynam:asymp}.
Indeed, it easily follows that the contribution to
\eqref{dynam:asymp} from {\em imprimitive} closed geodesics is
$$O\left(Te^{\frac{P(-\cH/2)T}{2}}\right)$$ and grows much slower than
$e^{P(-\cH/2)T}$. Here we have used that every imprimitive geodesic
corresponds to a primitive geodesic of at least twice smaller length.
Therefore, the asymptotics in \eqref{dynam:asymp} still
holds if we restrict summation to the primitive closed geodesics.
\qed

\end{document}